\documentclass[12pt,twoside]{article}
\usepackage{amsmath}
\usepackage{latexsym}
\usepackage{amssymb}
\usepackage{a4wide}
\newtheorem{thm}{Theorem}[section]
\newtheorem{la}[thm]{Lemma}
\newtheorem{Defn}[thm]{Definition}
\newtheorem{Remark}[thm]{Remark}
\newtheorem{Note}[thm]{Note}
\newtheorem{prop}[thm]{Proposition}

\newtheorem{Example}[thm]{Example}
\newtheorem{Examples}[thm]{Examples}
\newtheorem{Problems}[thm]{Problems}

\newtheorem{Problem}[thm]{Problem}
\newtheorem{Number}[thm]{\!\!}
\newenvironment{defn}{\begin{Defn}\rm}{\end{Defn}}

\newenvironment{example}{\begin{Example}\rm}{\end{Example}}

\newenvironment{rem}{\begin{Remark}\rm}{\end{Remark}}
\newenvironment{numba}{\begin{Number}\rm}{\end{Number}}
\newenvironment{proof}{{\noindent\bf Proof.}}%
                  {\nopagebreak\hspace*{\fill}$\Box$\medskip\medskip\par}   

\newcommand{\wb}{\overline}

\newcommand{\wt}{\widetilde}

\newcommand{\n}{\rm}

\newcommand{\mto}{\mapsto}
\newcommand{\emb}{\hookrightarrow}

\newcommand{\ve}{\varepsilon}
\newcommand{\isom}{\cong}

\newcommand{\N}{{\mathbb N}}

\newcommand{\R}{{\mathbb R}}
\newcommand{\bO}{{\mathbb O}}

\newcommand{\Z}{{\mathbb Z}}
\newcommand{\C}{{\mathbb C}}

\newcommand{\K}{{\mathbb K}}
\newcommand{\sys}{{\cal S}}

\newcommand{\F}{{\mathbb F}}
\newcommand{\A}{{\mathbb A}}

\newcommand{\cC}{{\cal C}}

\newcommand{\cP}{{\cal P}}
\newcommand{\cT}{{\cal T}}
\newcommand{\cW}{{\cal W}}
\newcommand{\g}{{\mathfrak g}}

\newcommand{\ch}{{\mathfrak h}}

\newcommand{\ck}{{\mathfrak k}}

\newcommand{\cs}{{\mathfrak s}}
\newcommand{\Cr}{{\mathfrak r}}
\newcommand{\dl}{{\displaystyle \lim_{\longrightarrow}}}

\newcommand{\Aut}{\mbox{\n Aut}}

\newcommand{\dsemi}{\mbox{$\times\!$\rule{.15 mm}{2.07 mm}\,}}

\newcommand{\sub}{\subseteq}

\newcommand{\GL}{\mbox{\rm GL}}

\newcommand{\gl}{{\mathfrak g}{\mathfrak l}}

\newcommand{\im}{\mbox{\n im}}

\newcommand{\car}{\mbox{{\rm char}}}

\newcommand{\diag}{\mbox{\rm diag}}

\newcommand{\id}{\mbox{\n id}}

\newcommand{\cB}{{\cal B}}
\newcommand{\cA}{{\cal A}}
\newcommand{\cD}{{\cal D}}

\newcommand{\cL}{{\cal L}}

\newcommand{\ob}{\mbox{\n ob}}

\newcommand{\smim}{\mbox{\n\footnotesize im}}

\newcommand{\Supp}{\mbox{\n supp}}

\newcommand{\sbull}{{\scriptscriptstyle \bullet}}

\newcommand{\evol}{{\mbox{{\rm evol}}}}
\newcommand{\Evol}{{\mbox{{\rm Evol}}}}

\newcommand{\Fix}{{\mbox{{\rm Fix}}}}
\newcommand{\ops}{{\mbox{{\rm\footnotesize op}}}}

\begin{document}
\begin{center}
{\Large \bf Fundamentals of Direct Limit Lie Theory}\vspace{4.5 mm}\\
{\bf Helge Gl\"{o}ckner}\vspace{.9mm}
\end{center}
\noindent{\bf Abstract.\/}
We show that
every countable direct system of
finite-dimensional real or complex
Lie groups
has a direct limit in the category
of Lie groups modelled on locally convex spaces.
This enables us to push all basic
constructions of finite-dimensional
Lie theory to the case of direct limit groups.
In particular, we obtain an analogue
of Lie's third theorem:
Every countable-dimensional
locally finite real or complex Lie algebra
arises as the Lie algebra
of some regular Lie group
(a suitable direct limit group).\\[3mm]
{\footnotesize
}
\section*{Introduction}
\noindent
In this paper, we develop the foundations of
Lie theory for countable
direct limits of finite-dimensional Lie groups.
For the purposes of this introduction,
consider an ascending sequence
$G_1\sub G_2\sub\cdots$ of finite-dimensional
real Lie groups, such that the inclusion
maps are smooth homomorphisms.
Then $G:=\bigcup_{n\in \N}G_n$
is a group in a natural way,
and it becomes a topological group when equipped
with the final topology with respect to the inclusion maps
$G_n\to G$ (\cite{HSTH}, \cite{TSH}).
A simple example is
$\GL_\infty(\R)$,
the group of invertible
matrices of countable size, differing from the unit
matrix at only finitely many places.
Here $\GL_\infty(\R)=\bigcup_n\GL_n(\R)$,
where $\GL_1(\R)\sub \GL_2(\R)\sub \cdots$
identifying $A\in \GL_n(\R)$ with $\diag(A,1)\in \GL_{n+1}(\R)$.
Our goal is
to make $G=\bigcup_n G_n$
a (usually infinite-dimensional) Lie group, and to discuss the
fundamental constructions of Lie theory
for such groups.\\[3mm]
{\bf Existing methods.}
Provided certain technical conditions
are satisfied
(ensuring in particular that
$\exp_G:=\dl\, \exp_{G_n}\!: \dl\,L(G_n)\to \dl\, G_n=G$\vspace{-.8mm}
is a local homeomorphism at~$0$),
the map $\exp_G$
restricts to a chart making $G$
a Lie group (see \cite{NRW1}, \cite{NRW2}
and \cite[Appendix]{NRW3}).
This method applies, in particular,
to $\GL_\infty(\R)$ and other direct limits
of linear Lie groups. It produces
Lie groups which are not only smooth,
but
real analytic in the sense
of convenient differential calculus \cite[Rem.\,6.5]{DIR}.
It is also known that every
Lie subalgebra of $\gl_\infty(\R):=\dl\, \gl_n(\R)$\vspace{-.8mm}
integrates to a subgroup of $\GL_\infty(\R)$
\cite[Thm.\,47.9]{KaM};
this provides an alternative construction
of the Lie group structure on
various direct limit groups.
But neither of these methods
is general enough to tackle
arbitrary direct limits of Lie groups.
In particular, examples show that
$\exp_G$ need not be injective
on any $0$-neighbourhood~\cite[Example~5.5]{DIR}.
Therefore a general construction
of a Lie group structure on $G=\bigcup_n G_n$
cannot make use of~$\exp_G$.

\noindent
{\bf A general construction principle.}
In \cite{DIR},
a smooth Lie group structure on $G=\bigcup_n G_n$
was constructed in the case where all
inclusion maps are embeddings (for
``strict'' direct systems).
Strict direct limits of Lie groups
are discussed there as special cases of direct limits
of direct sequences $M_1\sub M_2\sub\cdots$
of finite-dimensional smooth manifolds
and embeddings onto closed submanifolds.
To make $M:=\bigcup_n M_n$ a smooth manifold,
one starts with a chart $\phi_{n_0}$
of some $M_{n_0}$ and then uses tubular neighbourhoods
to extend $\phi_n$ already constructed
(possibly restricted to a smaller
open set) to a chart $\phi_{n+1}$ of $M_{n+1}$.
Then $\dl\, \phi_n$\vspace{-.8mm} is a chart for~$M$.
In the present article, we generalize
this construction principle
in two ways. First, we
are able to remove the strictness condition.
This facilitates to make $\bigcup_n M_n$
a smoothly paracompact,
smooth manifold, for any ascending sequence
of paracompact,
finite-dimensional smooth manifolds
and injective immersions
(Theorem~\ref{createmfd}, Proposition~\ref{mfde}).
Second, we generalize the method
from the case of smooth manifolds
over~$\R$ to the case of real- and complex
analytic manifolds
(Theorem~\ref{createmfd},
Proposition~\ref{mfdg}).
This enables us to turn $G:=\bigcup_n G_n$
into a
real analytic 
Lie group in the sense of convenient
differential calculus,
resp., a complex Lie group,
for any ascending sequence of finite-dimensional
real or complex Lie groups (Theorem~\ref{limgps}).\footnote{More
generally,
we can create direct limit Lie
groups for arbitrary countable direct
systems of finite-dimensional real or complex Lie groups.
The bonding maps need not be injective.}
Each direct limit group
$G$ is regular in the convenient sense
(the argument from \cite[Thm.\,47.8]{KaM} carries over).
Moreover, $G$ is
a regular
Lie group in Milnor's sense (Theorem~\ref{Milnorreg}):
this is much harder to prove.\\[3mm]
{\bf Lie theory for direct limit groups.}
Despite the fact that $\exp_G$ need not be well-behaved,
all of the basic constructions
of finite-dimensional Lie theory
can be pushed to the case of direct limit
groups $G=\bigcup_n G_n$.
Thus, subgroups and Hausdorff
quotient groups are Lie groups
(Propositions~\ref{subgplie}
and \ref{homog}),
a universal complexification $G_\C$
exists (Proposition~\ref{univcx}),
subalgebras of $L(G)$
integrate to analytic subgroups (Proposition~\ref{intsubalg}),\linebreak
and Lie algebra homomorphisms integrate to
group homomorphisms in the expected way (Proposition~\ref{inthom}).
Furthermore (Theorem~\ref{intthm}), every
locally finite real or complex
Lie algebra of countable dimension is enlargible,
i.e., it arises as the Lie algebra of some
Lie group
(a suitable direct limit group).
Such Lie algebras have been studied
by Yu. Bahturin, A.\,A. Baranov, I. Dimitrov, K.-H. Neeb,
I. Penkov, H. Strade, N. Stumme, A.\,E. Zalesskii,
and others.
If $H\sub G$ is a closed subgroup,
then $H$ is a
conveniently real analytic ($c^\omega_\R$-)
submanifold of~$G$.
Furthermore, the homogeneous space
$G/H$ can be given a $c^\omega_\R$-manifold
structure making $\pi\!: G\to G/H$ a $c^\omega_\R$-principal bundle
(Proposition~\ref{homog}).
Similar results are available
for complex Lie groups.
We remark that special cases of complexifications
and homogeneous spaces of direct limit groups
have already been used in \cite{NRW3} and \cite{Wol},
in the context of a Bott-Borel-Weil theorem,
resp., direct limits of principal series
representations. Universal complexifications
of ``linear'' direct limit groups $G\sub \GL_\infty(\R)$
have been discussed in~\cite{GCX},
in the framework of BCH-Lie groups.
For some special examples of direct limit
manifolds of relevance for
algebraic topology, see \cite[\S47]{KaM}.\\[3mm]
{\bf Variants.}
Although our main results concern the real
and complex cases,
some of the constructions apply just as well
to Lie groups over local fields
(i.e., totally disconnected,
locally compact, non-discrete topological fields,
such as the $p$-adic numbers), and are formulated accordingly.
Readers mainly interested in
the real and complex cases are invited
to read ``$\K$'' as $\R$ or $\C$,
ignore the definition of smooth maps
over general topological fields,
and assume that all Lie groups are modelled
on real or complex locally convex spaces.\vspace{-2mm}
\section{Basic definitions and facts}
We are working
in two settings of differential calculus
in parallel: 1.\ The Convenient Differential Calculus
of Fr\"{o}licher, Kriegl and Michor.
2.\ Keller's $C^\infty_c$-theory
(going back to Michal and Bastiani),
as used, e.g., in \cite{MiP}, \cite{Mil},
\cite{GCX}, \cite{RES},
and generalized to a general differential calculus
over topological fields in~\cite{BGN}.
For the basic notions of infinite-dimensional Lie theory
($L(G)$, $\exp_G$, logarithmic derivative,
product integral), see
\cite{KaM} and \cite{Mil}.
\begin{numba}
{\bf Convenient differential calculus.}
Our source
for Convenient Differential Calculus is \cite{KaM},
and we presume familiarity with the basic ideas.
The smooth maps and manifolds
from convenient calculus will be called
{\em $c^\infty_\R$-maps\/}
and {\em $c^\infty_\R$-manifolds\/}
here. Maps and manifolds
which are holomorphic in the convenient sense
will be called $c^\infty_\C$ or $c^\omega_\C$.
Real analytic maps and manifolds in the convenient
sense will be called $c^\omega_\R$.
Likewise for Lie groups. The
regular $c^\infty_\R$-Lie groups from convenient
calculus (see \cite[Defn.\,38.4]{KaM})
will be called {\em $c^\infty_\R$-regular\/};
we call
a $c^\infty_\C$-Lie group
{\em $c^\infty_\C$-regular\/} or {\em $c^\omega_\C$-regular\/}
if its underlying
$c^\infty_\R$-Lie group is $c^\infty_\R$-regular.
A $c^\omega_\R$-Lie group $G$ will be called
{\em $c^\omega_\R$-regular\/} if it is $c^\infty_\R$-regular
and the right product integral
$\Evol^r_G(\gamma)\!: \R\to G$
of each $c^\omega_\R$-curve $\gamma\!: \R\to L(G)$
is~$c^\omega_\R$.
The definitions of $c^\omega_\R$-regularity
and $c^\omega_\C$-regularity ensure the following:
\end{numba}
\begin{la}\label{useofreg}
Given $\K\in \{\R,¸\C\}$,
suppose that
$G$ and $H$ are $c^\omega_\K$-Lie groups,
where $G$ is simply connected and $H$ is
$c^\omega_\K$-regular.
Then, for every bounded $\K$-Lie algebra
homomorphism $\alpha\!: L(G)\to L(H)$,
there exists a unique $c^\omega_\K$-homomorphism
$\beta\!: G\to H$ such that $L(\beta)=\alpha$.
\end{la}
\begin{proof}
By \cite[Thm.\,40.3]{KaM},
there exists a unique $c^\infty_\R$-homomorphism
$\beta\!: G\to H$ such that $L(\beta)=\alpha$.
If $\K=\R$ and $\gamma\!: \R\to G$ is a $c^\omega_\R$-curve,
then $\beta\circ \gamma\!: \R\to H$ is a smooth curve
with right logarithmic derivative
$\delta^r(\beta\circ \gamma)=L(\beta)\circ \delta^r\gamma=\alpha\circ
\delta^r\gamma$.
Here $\alpha\circ \delta^r\gamma$ is~$c^\omega_\R$,
whence its right product integral $\beta\circ \gamma$
is~$c^\omega_\R$, by
$c^\omega_\R$-regularity. Hence $\beta$ is $c^\omega_\R$.
If $\K=\C$, then $\beta$ is a $c^\infty_\R$-homomorphism
such that $T_x(\beta)$ is $\C$-linear for each $x\in G$,
as $T_1(\beta)=\alpha$ is $\C$-linear.
Hence $\beta$ is $c^\omega_\C$ by \cite[Thm.\,7.19\,(8)]{KaM}.\vspace{-3.5mm}
\end{proof}
\begin{numba}{\bf Keller's {\boldmath $C^\infty_c$\/}-theory
and analytic maps.}\label{Keller}
Let $E$ and $F$ be locally convex spaces
over $\K\in \{\R,\C\}$,
$U\sub E$ be open and
$f\!: U\to F$ be a map.
If $\K=\R$ and $r\in \N_0\cup\{\infty\}$,
then $f$ is called $C^r_\R$ if it is continuous
and,
for all $k\in \N_0$ such that $k\leq r$,
the iterated directional derivatives $d^kf(x,y_1,\ldots,y_k):=
D_{y_1}\cdots D_{y_k}f(x)$ exist
for all $x\in U$ and $y_1,\ldots,y_k\in E$,
and define a continuous map $d^kf\!: U\times E^k\to F$.
The $C^\infty_\R$-maps are also called {\em smooth}.
If $\K=\C$, we call $f$ a {\em $C^\infty_\C$-map},
$C^\omega_\C$, or {\em complex analytic\/},
if it is continuous and given locally by a
pointwise convergent series of continuous
homogeneous polynomials~\cite[Defn.\,5.6]{BaS}.
If $\K=\R$,  we call $f$ {\em real analytic\/} or $C^\omega_\R$
if it extends to a complex analytic map between
open subsets of the complexifications of~$E$ and~$F$.
See \cite{MiP}, \cite{Mil}, or
\cite{RES}
for further information
(also concerning the corresponding smooth
and $\K$-analytic
Lie groups and manifolds).
\end{numba}
\begin{numba}{\bf General differential calculus.}
Let $E$ and $F$ be (Hausdorff) topological vector spaces
over a non-discrete topological field~$\K$,
$U\sub E$ be open, and $f\!: U\to F$ a map.
According to \cite{BGN},
$f$ is called $C^1_\K$
if it is continuous and there
exists a (necessarily unique)
continuous map $f^{[1]}\!: U^{[1]}\to F$
on $U^{[1]}:=\{(x,y,t)\in U\times E\times \K\!:
x+ty\in U\}$ such that $f^{[1]}(x,y,t)=\frac{1}{t}(f(x+ty)-f(x))$
for all $(x,y,t)\in U^{[1]}$ such that $t\not=0$.
Inductively, $f$ is called $C^{k+1}_\K$ if it is
$C^1_\K$ and $f^{[1]}$ is $C^k_\K$;
it is  $C^\infty_\K$ if it is $C^k_\K$
for all~$k$.
As shown in \cite{BGN},
compositions of $C^k_\K$-maps
are $C^k_\K$, and
being $C^k_\K$ is a local property.
For
maps between open subsets
of locally convex spaces,
the present definitions
of $C^k_\R$-maps and $C^\infty_\C$-maps
are equivalent to those from {\bf \ref{Keller}}
(\cite{BGN}, Prop.\ 7.4 and 7.7).
Analytic maps between
open subsets of Banach spaces over a complete
valued field $\K$ (as used in \cite{Bou} or \cite{Ser})
are $C^\infty_\K$ \cite[Prop.\,7.20]{BGN}.
For further
information, also
concerning $C^\infty_\K$-manifolds
and Lie groups modelled on topological $\K$-vector spaces,
we refer to \cite{BGN}, \cite{IMP},
\cite{ANA}, and \cite{GEN}.
\end{numba}
\begin{numba} {\bf Direct limits.}
A {\em direct system\/}
in a category $\A$ is a pair
$\sys=((X_i)_{i\in I},(\phi_{i,j})_{i\geq j})$,
where $(I,\leq)$ is a directed set,
each $X_i$ an object of~$\A$,
and each $\phi_{i,j}\!: X_j\to X_i$ a morphism
(``bonding map'')
such that $\phi_{i,i}=\id_{X_i}$
and $\phi_{i,j}\circ \phi_{j,k}=\phi_{i,k}$ if
$i\geq j\geq k$. A {\em cone over $\sys$\/}
is a pair $(X,(\phi_i)_{i\in I})$, where $X\in \ob \,\A$
and $\phi_i\!: X_i\to X$ is a morphism
for $i\in I$ such that $\phi_i\circ \phi_{i,j}=\phi_j$
if $i\geq j$. A cone $(X,(\phi_i)_{i\in I})$
is a {\em direct limit cone\/} over $\sys$
in the category~$\A$
if, for every cone $(Y,(\psi_i)_{i\in I})$ over~$\sys$,
there exists a unique morphism $\psi\!: X\to Y$
such that $\psi\circ \phi_i=\psi_i$ for each~$i$.
We then write $(X,(\phi_i)_{i\in I})=\dl\, \sys$.\vspace{-.8mm}
If the bonding maps and ``limit maps''
$\phi_i$ are understood,
we simply call $X$ the {\em direct limit\/} of $\sys$
and write $X=\dl\, X_i$.\vspace{-.8mm}
If also $\cT=((Y_i)_{i\in I},(\psi_{i,j})_{i\leq j})$
is a direct system over $I$ and
$(Y,(\psi_i)_{i\in I})$ a cone over $\cT$,
we call a family $(\eta_i)_{i\in I}$
of morphisms $\eta_i\!: X_i\to Y_i$
{\em compatible\/}
if $\eta_i\circ \phi_{i,j}=\psi_{i,j}\circ \eta_j$
for $i\geq j$.
Then $(Y,(\psi_i\circ \eta_i)_{i\in I})$
is a cone over $\sys$;
write $\dl\,\eta_i:=\eta$\vspace{-.8mm}
for the morphism $\eta\!: X\to Y$
such that $\eta\circ \phi_i=\psi_i\circ \eta_i$.
If there is a compatible family
$(\eta_i)_{i\in I}$ with each $\eta_i$ an isomorphism,
$\sys$ and $\cT$
are called {\em equivalent}.
Then $\sys$ has a direct limit if
and only if so does $\cT$; in this case,
$\dl\,\eta_i$\vspace{-.8mm}
is an isomorphism. Every countable direct
set has a {\em cofinal subsequence\/},
whence countable direct systems
can be replaced by {\em direct sequences},
viz.\ $I=(\N,\leq)$.
\end{numba}
\begin{numba}\label{dirset}
{\bf Direct limits of sets, topological spaces, and groups.}
If $\sys=((X_i)_{i\in I},(\phi_{i,j})_{i\geq j})$
is a direct system of sets,
write $(j,x)\sim (k,y)$ if there exists $i\geq j,k$
such that $\phi_{i,j}(x)=\phi_{i,k}(y)$;
then $X:=\big(\coprod_{i\in I}X_i\big)/\!\sim$, together with the
maps $\phi_i\!: X_i\to X$, $\phi_i(x):=[(i,x)]$,
is the direct limit of $\sys$ in the category of sets.
Here $X=\bigcup_{i\in I}\phi_i(X_i)$.
If each $\phi_{i,j}$ is injective, then so is
each $\phi_i$, whence $\sys$ is equivalent to
the direct system of the subsets $\phi_i(X_i)\sub X$,
together with the inclusion maps.
This facilitates to replace injective direct
systems by direct systems in which all bonding maps
are inclusion maps.
If $\sys:=((X_i)_{i\in I},(\phi_{i,j}))$
is a direct
system of topological spaces and continuous maps,
then the direct limit $(X,(\phi_i)_{i\in I})$ of the underlying sets
becomes the direct limit in the category of topological spaces
and continuous maps if we equip $X$ with the {\em
DL-topology}, the final topology with respect to the family
$(\phi_i)_{i\in I}$. Thus $U\sub X$ is open if and
only if $\phi^{-1}(U)$ is open in $X_i$, for each~$i$.
If $\sys$ is {\em strict\/} in the sense that each $\phi_{i,j}$
is a topological embedding,
then also each $\phi_i$ is a topological embedding \cite[La.\,A.5]{NRW2}.
If $((G_i)_{i\in I},(\phi_{i,j})_{i\geq j})$
is a direct system of groups and homomorphisms,
then the
direct limit $(G,(\phi_i)_{i\in I})$ of the underlying sets
becomes the direct limit in the category of
groups and homomorphisms when equipped with the unique
group structure making each $\phi_i$ a homomorphism;
the group inversion and multiplication on $G$
are $\dl\, \kappa_i$\vspace{-3mm} and $\dl\,\mu_i$,
in terms of those on the $G_i$'s.
\end{numba}
For further information concerning direct limits
of topological groups and topological spaces,
see \cite{DIR}, \cite{Han}, \cite{HSTH}, and
\cite{TSH}.
\begin{la}\label{base}
Let $X_1\sub X_2\sub \cdots$ be an
ascending sequence of topological spaces such that the inclusion maps
are continuous; equip $X:=\bigcup_{n\in \N}X_n$
with the final topology with respect to the
inclusion maps $\lambda_n\!: X_n\to X$
$($the DL-topology$)$.
Then the following holds:
\begin{itemize}
\item[\rm (a)]
If each $X_n$ is $T_1$,
then so is $X$.
\item[\rm (b)]
If $U_n\sub X_n$ is open and $U_1\sub U_2\sub \cdots$,
then $U:=\bigcup_n U_n$ is open in~$X$ and the
DL-topology on $U=\dl\,U_n$\vspace{-1mm}
coincides with the topology induced by~$X$.
\item[\rm (c)]
If each $X_n$ is locally compact, then $X$ is Hausdorff.
\item[\rm (d)]
If each $X_n$ is $T_1$ and $K\sub X$ is compact, then
$K\sub X_n$ for some~$n$.
\end{itemize}
\end{la}
\begin{proof}
(a) Let $x\in X$. Then $\lambda_n^{-1}(\{x\})$ is either
$\{x\}$ or empty, hence closed in the $T_1$-space~$X_n$.
Hence $\{x\}$ is closed in $X$.\vspace{1mm}

(b) and (c): This is proved in
\cite[Prop.\,4.1\,(ii)]{Han}
and \cite[La.\,3.1]{DIR}
for strict direct sequences,
but the strictness is not used in the proofs.\vspace{1mm}

(d) If not, for each $n$ we find $x_n\in K\setminus X_n$.
Then $D:=\{x_n\!: n\in \N\}\sub K$
is closed in~$X$ (and thus compact),
as $D\cap X_n$ is finite and thus
closed,
for each~$n$. On the other hand,
$D=\dl\, (D\cap X_n)$\vspace{-.8mm}
for the topology induced by~$X$,
as $D$ is closed in~$X$.
Now $D\cap X_n$ being discrete, this entails
$D$ is discrete and hence finite
(being also compact). Contradiction.\vspace{-1mm}
\end{proof}
\begin{numba}\label{new1.8}
Let $E$ be a countable-dimensional
vector space over a
non-discrete, locally compact topological field~$\K$
(e.g., $\K=\R$ or~$\C$).
Then the finest vector topology on~$E$
is locally convex and coincides
with the so-called {\em finite topology},
the final topology with respect to
the inclusion maps $F\to E$,
where $F$ ranges through
the set of finite-dimensional
vector subspaces of~$E$
(and $F$ is
equipped with its
canonical Hausdorff vector topology).
Thus, the finite topology on~$E$ is
the DL-topology on
$E=\dl\,F$.\vspace{-.8mm}
See \cite{DIR} and the references therein for these standard
facts.
The space $\K^\infty:=\K^{(\N)}=\dl\, \K^n$\vspace{-.8mm}
of finite sequences will always be
equipped with the finite topology.
We shall frequently
identify
$\K^n$ with the subspace
$\K^n\times \{0\}$ of $\K^\infty$,
and $\K^m$ with $\K^m\times\{0\}\sub \K^n$ if $n\geq m$. 
\end{numba}
\begin{la}\label{silva}
Let $\K$ be $\R$, $\C$ or a local field,
and $E$ be a $\K$-vector space
of countable dimension, equipped
with the finite topology.
Let $E_1\sub E_2\sub \cdots$
be an ascending sequence of vector subspaces of~$E$
such that $\bigcup_{n\in \N}E_n=E$,
and $U_n\sub E_n$ be open subsets
such that $U_1\sub U_2\sub \ldots$.
Let $f\!: U\to F$
be a map into a topological $\K$-vector space~$F$
on the open subset
$U:=\bigcup_{n\in \N}U_n$ of $E$.
Then the following holds:
\begin{itemize}
\item[\rm (a)]
Given $r\in \N_0\cup\{\infty\}$,
$f$ is $C^r_\K$
if and only if $f_n:=f|_{U_n}\!: U_n\to F$ is $C^r_\K$
for each~$n$.
\item[\rm (b)]
If $\K\in \{\R,\C\}$
and $F$ is locally convex and Mackey complete,
then $f$ is $C^\infty_\K$
if and only if it is $c^\infty_\K$.
Furthermore,
$f$ is $c^\omega_\K$
if and only if $f|_{U_n}$
is $c^\omega_\K$ for each $n\in \N$.
\end{itemize}
\end{la}
\begin{proof}
(a)\,\footnote{For $\K=\R$ and
locally convex~$F$,
see also \cite{DIR}, lines preceding
La.\,4.1. This implies the claim
for $\K=\C$, $r=\infty$, $F$ locally convex
because then $f$ is $C^\infty_\R$
with $df(x,\sbull)$
complex linear for each~$x$
(because $df(x,\sbull)|_{E_n}=df_n(x,\sbull)$), whence $f$ is complex
analytic by~\cite[La.\,2.5]{RES}.}
We may assume $r<\infty$.
Lemma~\ref{base}\,(b)
settles the case $r=0$.
If $r\geq 1$, note that
$U_1^{[1]}\sub U^{[1]}_2\sub\cdots$
and $U^{[1]}=\bigcup_n U_n^{[1]}$.
The product topology on
$E\times E\times \K$
is the finite topology
(cf.\ \cite[Prop.\,3.3]{DIR})
and hence induces on $U^{[1]}$
the topology making it the direct
limit topological space $U^{[1]}=\dl\, U^{[1]}_n$\vspace{-1.8mm}
(Lemma~\ref{base}\,(b)).
By induction, the cone $(F,(f_n^{[1]})_{n\in \N})$
of $C^{r-1}_\K$-maps induces a
$C^{r-1}_\K$-map $g\!: U^{[1]}=
\dl\, U_n^{[1]}\to F$,\vspace{-1.8mm}
determined by $g|_{U_n^{[1]}}=f_n^{[1]}$.
As $g$ is continuous
and extends the
difference quotient map,
$f$ is $C^1_\K$ with $f^{[1]}=g$.\vspace{-1.5mm}
Now $f$ being
$C^1_\K$ with $f^{[1]}=g$ of class $C^{r-1}_\K$,
the map~$f$ is $C^r_\K$.\vspace{1mm}

(b) If $f$ is $C^\infty_\K$, then
it is $c^\infty_\K$.
If $f$ is $c^\infty_\K$, then $f|_{U_n}$ is $c^\infty_\K$
for each $n$ and thus $C^\infty_\K$, as
$\dim_\K(E_n)<\infty$.
Hence $f$ is $C^\infty_\K$,
by (a).
Given a $c^\omega_\R$-curve $\gamma\!: \R\to U$
and $t_0\in \R$, pick an open
relatively compact neighbourhood $J\sub \R$ of~$t_0$.
Then $\gamma(J)\sub U_n$ for some
$n$ by Lemma~\ref{base}\,(d),
and thus $\gamma|_J$ is $c^\omega_\R$
if so is~$f|_{U_n}$. The remainder is now obvious.\vspace{-.5mm}
\end{proof}
A map $f\!: \R^\infty\to \R$
which is $C^\omega_\R$ on each $\R^n$
need not be $C^\omega_\R$ \cite[Ex.\,10.8]{KaM}.
For this reason, we have to work
with the weaker concept of $c^\omega_\R$-maps.
\section{Extension of charts}
In this section, we explain how
a chart of a submanifold $M_1\sub M_2$
(or its restriction to a slightly
smaller open set)
can be extended to a chart
of~$M_2$.
\begin{la}\label{smoothext}
Let $M_1$ and $M_2$ be finite-dimensional smooth
$($resp., real analytic$)$
manifolds over~$\R$,
of dimensions $m_1$ and $m_2$, respectively.
Assume that $M_1\sub M_2$
and assume that the inclusion map $\lambda \!: M_1\to M_2$
is a smooth $($resp.,
real analytic$)$ immersion. Let $\phi_1\!: U_1 \to V_1$
be a chart of~$M_1$, where $U_1$ is open in $\R^{m_1}$
and $V_1$ is an open, relatively compact,
contractible subset of~$M_1$.
Then there exists a chart $\phi_2\!: U_2\to V_2$
of $M_2$ such that $U_2\cap (\R^{m_1}\times \{0\}) =U_1\times\{0\}$,
$\phi_2(x,0)=\phi_1(x)$ for all $x\in U_1$,
and such that $V_2\sub M_2$ is relatively compact
and contractible.
\end{la}
\begin{proof}
Because $C:=\wb{V_1}\sub M_1$ is compact,
the map $\lambda|_C$ is a topological embedding.
Now $V_1$ being open in~$C$,
we deduce that $V_1=\lambda|_C(V_1)$
is open in $\lambda(C)$,
whence there exists an open subset $W\sub M_2$
such that $W\cap \lambda(C)=V_1$.
Since $\lambda(C)$ is closed in $M_2$,
the preceding formula shows that $V_1$ is closed in~$W$.
After shrinking $W$, we may assume that $W$ is
$\sigma$-compact, and relatively compact in~$M_2$.
Then $V_1$ is a closed submanifold
of the $\sigma$-compact, relatively compact,
open submanifold $W$ of~$M_2$.
{\em Smooth case\/}:
By \cite[Thm.\,IV.5.1]{Lan},
$V_1$ admits a smooth tubular neighbourhood
in~$W$, i.e.,
there exists a $C^\infty_\R$-diffeomorphism
$\psi\!: V_2\to P$
from some open neighbourhood $V_2$ of $V_1$
in~$W$ onto some open neighbourhood $P$
of the zero-section of some smooth
vector bundle $\pi\!: E\to V_1$
over $V_1$, such that $\psi|_{V_1}=\id_{V_1}$
(identifying $V_1$ with the zero-section of~$E$).
{\em Real analytic case\/}:
Being $\sigma$-compact,
$W$ is $C^\omega_\R$-diffeomorphic
to a closed real analytic submanifold of $\R^k$
for some~$k\in \N_0$ (see \cite[Thm.\,3]{Gr2}),
whence $W$ admits a real analytic Riemannian
metric~$g$.
Using the real analytic
Riemannian metric, the classical construction
of tubular neighbourhoods
provides a real analytic
tubular neighbourhood $\psi\!: V_1\supseteq
V_2\to P\sub E$.

In either case,
after shrinking $V_2$ and $P$, we may assume
that $P$ is balanced, i.e.,
$[-1,1]P\sub P$ (using the scalar multiplication
in the fibres of~$E$).
Being a vector bundle over a contractible,
$\sigma$-compact base manifold, $E$ is trivial.
This is well-known
in the smooth case \cite[Cor.\,4.2.5]{Hir}.
For the real analytic case, note that
$E$ is associated to a real analytic $\GL(F)$-principal
bundle over the $\sigma$-compact, contractible
$C^\omega_\R$-manifold~$V_1$, where
$F:=\R^{m_2-m_1}$ is the fibre of~$E$.
This principal bundle is trivial by \cite[Teorema~5]{Tog}
(combined with \cite[Cor.\,4.2.5]{Hir}),
and hence so is~$E$.
(Compare also \cite{Anc} and \cite{Gua}).

By the preceding, we
find an isomorphism of smooth (resp., real
analytic) vector bundles
$\theta\!: E\to V_1\times \R^{m_2-m_1}$.
Then $\kappa\!: \phi_1^{-1}\times \id\!: V_1\times\R^{m_2-m_1}
\to U_1\times \R^{m_2-m_1}\sub \R^{m_2}$
is a $C^\infty_\R$- (resp., $C^\omega_\R$-) diffeomorphism,
and $U_2:=\kappa(\theta(P))$ is an open subset of $\R^{m_2}$
such that $U_2\,\cap \,(\R^{m_1}\!\times\! \{0\})=U_1$.
Then $\phi_2:=(\kappa\circ \theta\circ \psi)^{-1}|_{U_2}^{V_2}\!:
U_2\to V_2$ is a $C^\infty_\R$- (resp., $C^\omega_\R$-) diffeomorphism
from $U_2$ onto the open subset $V_2$ of~$M_2$,
such that $\phi_2(x,0)=\phi_1(x)$ for all $x\in U_1$.
Since $V_2\sub W$, the set $V_2$ is relatively compact in~$M_2$.
To see that $V_2$ is contractible,
we only need to show that so is~$P$,
as $V_2$ and $P$ are homeomorphic.
Let $H\!: [0,1]\times V_1\to V_1$
be a homotopy from $\id_{V_1}$ to a constant map.
The map $[0,1]\times P\to P$,
$(t,x)\mto (1-t)x$ (which uses scalar multiplication
in the fibres) is a homotopy from $\id_P$
to $\pi|_P$.
The map $[0,1]\times P\to P$,
$(t,x)\mto H(t,\pi(x))$ is a homotopy from
$\pi|_P$ to a constant map.
Thus $\id_P$ is homotopic to a constant map
and thus $P$ is contractible.\vspace{-2mm}
\end{proof}
\begin{defn}
Let $\K$ be $\R$, $\C$ or a local field,
and $|.|$ be an absolute value on~$\K$ defining its
topology.
Given $n\in \N_0$ and $r>0$, we let
\[
\Delta^n_r:=\{(x_1,\ldots,x_n)\in \K^n\!: \,\mbox{$|x_j|<r$
for all $j=1,\ldots, n$}\,\}
\]
be the $n$-dimensional polydisk
of radius $r$ around~$0$.
If we wish to emphasize the ground field,
we also write $\Delta^n_r(\K)$ for $\Delta^n_r$.
\end{defn}
If $\K$ is a local field,
we define
$C^\infty_\K$-immersions
(and $C^\infty_\K$-submersions)
between finite-dimensional $C^\infty_\K$-manifolds
analogous to the $\K$-analytic case~\cite{Ser}.
Because an
Inverse Function Theorem
holds for $C^\infty_\K$-maps~\cite{IMP},
$C^\infty_\K$-immersions and submersions
have the usual properties.
\begin{la}[Extension Lemma]\label{extpoly}
Let $\K$ be $\R$, $\C$ or a local field.
Let $M$ be a finite-dimensional
$C^\infty_\K$-manifold $($or a finite-dimensional
real analytic manifold$\,)$,
of dimension $m\in \N_0$,
and $\phi\!: \Delta^n_r\to M$ be a
$C^\infty_\K$ $($resp., real analytic$)$
injective immersion, where $n\in \{0,1,\ldots, m\}$
and $r\!>\!0$.
Then, for every $s\in \,]0,r[\,$,
there exists a $C^\infty_\K$-diffeomorphism
$($resp., a real analytic diffeo\-morphism$)$
$\psi\!: \Delta^m_s\to V$ onto an open
subset $V$ of~$M$
such that $\psi(x,0)=\phi(x)$ for all $x\in \Delta^n_s$.
If $\K$ is a local field, the conclusion remains
valid for $s=r$.
The subset $V\sub M$ can be chosen relatively compact.
\end{la}
\begin{proof}
Let $s\in \;]0,r[$ and $t\in \;]s,r[$.
\vspace{1mm}

{\em The case of smooth or analytic manifolds over $\K=\R$.}
We equip $M_1:=\phi(\Delta^n_r)$
with the smooth (resp., real analytic) manifold structure making
$\phi|^{M_1}\!: \Delta^n_r\to M_1$ a diffeomorphism.
Then the inclusion map $\lambda\!: M_1\to M$
is an immersion,
$V_1:=\phi(\Delta^n_t)$ is a relatively compact,
contractible, $\sigma$-compact open subset of $M_1$,
and $\phi_1:=\phi|_{\Delta_t^n}^{V_1}\!: \Delta^n_t\to V_1$
is a chart for $M_1$.
By Lemma~\ref{smoothext},
there exists a $C^\infty_\R$- (resp., $C^\omega_\R$-)
diffeomorphism
$\phi_2\!: U_2\to V_2$ from an open subset $U_2$ of $\R^m$
onto an open subset $V_2$ of~$M$ such that
$U_2\cap (\R^n\!\times\! \{0\})=\Delta^n_t\!\times\! \{0\}$
and $\phi_2(x,0)=\phi_1(x)=\phi(x)$
for all $x\in \Delta^n_t$.
Now $\wb{\Delta^n_s}\sub \R^n$ being compact, we find
$\ve>0$ such that $\wb{\Delta^n_s}\times \Delta^{m-n}_\ve\sub
U_2$. Then
\[
\psi\!: \Delta_s^m\to M, \quad \psi(x,y):=
\phi_2(x,{\textstyle \frac{\ve}{s}}y)\qquad
\mbox{for $x\in \Delta^n_s$, $y\in \Delta^{m-n}_s$}
\]
is a mapping with the required properties.\vspace{1mm}

{\em The case $\K=\C$.}
The map $\phi|_{\Delta^n_t}$ is an embedding of complex
manifolds, and hence so is $f\!: \Delta^n_1\to M$,
$f(x):=\phi(tx)$. By \cite[Prop.\,1]{Roy},
there exists a holomorphic embedding
$F\!: \Delta^m_{s/t}\times \Delta^{n-m}_1\to M$
such that $F(x,0)=f(x)$ for all $x\in \Delta^n_{s/t}$.
Then $\psi\!: \Delta^m_s\to M$, $\psi(x,y):=F(\frac{1}{t}x,\frac{1}{s}y)$
(where $x\in \Delta^n_s$, $y\in \Delta^{m-n}_s$)
is a holomorphic embedding with the desired properties.\vspace{1mm}

{\em Relative compactness of $V$.}
By the real or complex case
already discussed,
there exists an extension $\wt{\psi}\!: \Delta^m_t\to \wt{V}$
of $\phi|_{\Delta^n_t}$. Then $V:=\wt{\psi}(\Delta^m_s)$
is a relatively compact open subset of~$M$,
and $\psi:=\wt{\psi}|_{\Delta^m_s}^V$
has the desired properties.\vspace{1mm}

{\em The case where $\K$ is a local field.}
In this case, $\Delta^n_r$ is compact,
whence $\phi$ is a $C^\infty_\K$-diffeomorphism
from $\Delta^n_r$ onto the compact
$C^\infty_\K$-submanifold $M_1:=\im\, \phi$ of~$M$.
The proof of \cite[La.\,8.1]{DIR} (tackling the
$\K$-analytic case) carries over verbatim
to the case of $C^\infty_\K$-manifolds;
we therefore find a $C^\infty_\K$-diffeomorphism
$\theta\!: \Delta^n_r\times \bO^{m-n}\to M$
such that $\theta(x,0)=\phi(x)$,
where $\bO$ is the maximal compact subring of~$\K$.
Pick $a\in \K^\times$ such that $a \Delta^{m-n}_r\sub \bO^{m-n}$;
then $\psi\!: \Delta^m_r\to M$, $\psi(x,y):=\theta(x,ay)$
for $x\in \Delta^n_r$, $y\in \Delta^{m-n}_r$
(resp., its restriction to $\Delta^m_s$)
is the required chart for~$M$.
\end{proof}
\section{Direct limits of finite-dimensional manifolds}
Let $\K$ be $\R$, $\C$ or a local field.
Throughout this section,
we let $\sys:=
((M_i)_{i\in I},(\lambda_{i,j})_{i\geq j})$
be a direct system of finite-dimensional
$C^\infty_\K$-manifolds $M_i$
and injective $C^\infty_\K$-immersions $\lambda_{i,j}\!: M_j\to M_i$.
We let $(M,(\lambda_i)_{i\in I})$
be the direct limit of~$\sys$
in the category of topological spaces,
and abbreviate $s:=\sup\{ \dim_\K (M_i) \!: i\in I\}\in \N_0\cup\{\infty\}$.
Our goal is to make $M$ a manifold,
and study its basic properties.
\begin{thm}\label{createmfd}
There exists a uniquely determined
$C^\infty_\K$-manifold structure
on $M$, modelled on the complete, locally convex
topological $\K$-vector space $\K^s$,
which makes $\lambda_i\!: M_i\to M$ a $C^\infty_\K$-map,
for each $i\in I$, and such that $(M,(\lambda_i)_{i\in I})=
\dl\, \sys$\vspace{-.8mm} in the category of $C^\infty_\K$-manifolds
modelled on topological $\K$-vector spaces
$($and $C^\infty_\K$-maps$)$.
For each $i\in I$
and $x\in M_i$,
the differential $T_x(\lambda_i)\!: T_x(M_i)\to
T_{\lambda_i(x)}(M)$ is injective.
For each $r\in \N_0$,
the $C^r_\K$-manifold underlying~$M$ satisfies
$(M,(\lambda_i)_{i\in I})=
\dl\, \sys$\vspace{-.8mm} in the category of $C^r_\K$-manifolds
modelled on topological $\K$-vector spaces.
\end{thm}
\begin{proof}
After passing to a cofinal subsequence
of an equivalent direct system
(cf.\ {\bf \ref{dirset}}),
we may assume
without loss of generality that
$I=\N$,
$M_1\sub M_2\sub \cdots$,
and that the immersion $\lambda_{n,m}$ is the inclusion map
for all $n,m\in \N$
such that $n\geq m$.
We let $M:=\bigcup_{n\in \N} M_n$,
equipped with the final topology with respect to the inclusion maps
$\lambda_n\!: M_n\to M$; then $(M,(\lambda_n)_{n\in \N})=
\dl\, ((M_n),(\lambda_{n,m}))$\vspace{-.8mm} in the category of topological
spaces.
We abbreviate $d_n:=\dim_\K (M_n)$
and $c_n:=d_{n+1}-d_n$.\vspace{1mm}

Let $\cA$ be the set of all
maps $\phi\!: P_\phi\to Q_\phi\sub M$
such that $P_\phi=\bigcup_{n\in \N}P_n\sub \K^s$,
$Q_\phi=\bigcup_{n\in \N}Q_n$, and
$\phi=\dl\, \phi_n$\vspace{-.8mm}
for some sequence
$(\phi_n)_{n\in \N}$ of charts $\phi_n\!: P_n\to Q_n\sub M_n$,
where each $P_n$ is an open (possibly empty)
subset of $\K^{d_n}$,
$Q_n$ open in $M_n$,
and $Q_m\sub Q_n$
and $\phi_n|_{Q_m}=\phi_m$
whenever $n\geq m$.
Here Lemma~\ref{base}\,(b) allows us to
interpret the open subsets $P_\phi\sub \K^s$
and $Q_\phi\sub M$
as the direct limits $\dl\, Q_n$\vspace{-.8mm} and $\dl\, P_n$
in the category of topological spaces,
whence $\phi$ is continuous.
Because each $\phi_n$ is injective,
also $\phi$ is injective,
and furthermore $\phi$
is surjective,
by definition of~$Q_\phi$.
If $V\sub P_\phi$ is open,
then $V\cap P_n$ is open in $P_n$,
whence $S_n:=\phi_n(V\cap P_n)$
is open in $Q_n$.
Because $S_1\sub S_2\sub \cdots$,
the union $\phi(V)=\bigcup_{n\in \N}S_n$
is open in $Q_\phi $ (Lemma~\ref{base}\,(b)).
Thus $\phi$ is an open map.
We have shown that $\phi$ is a homeomorphism.\vspace{1mm}

We claim that $\cA$ is a $C^\infty_\K$-atlas
for~$M$.
We first show that $\bigcup_{\phi\in \cA}
Q_\phi=M$.
To this end, let $x\in M$.
Then there exists $\ell\in \N_0$ such that $x\in M_\ell$.
Define $r_n:=1+2^{-n}$ for $n\in \N$.
We let $\phi_n\!: P_n\to Q_n$ be the chart of~$M_n$
with $P_n:=Q_n:=\emptyset$, for all $n<\ell$.
We pick a chart $\psi_\ell\!: \Delta^{d_\ell}_{r_\ell}(\K)\to W_\ell\sub
M_\ell$ of $M_\ell$ around~$x$, such that $\psi_\ell(0)=x$.
Inductively, the Extension Lemma~\ref{extpoly}
provides charts $\psi_n\!: \Delta^{d_n}_{r_n}\to W_n\sub M_n$
for $n\in \{\ell+1,\ell+2,\ldots\}$
such that $\psi_n|_{\Delta^{d_{n-1}}_{r_n}}=
\psi_{n-1}|_{\Delta^{d_{n-1}}_{r_n}}$\vspace{-.8mm} (identifying
$\K^{d_{n-1}}$ with $\K^{d_{n-1}}\times \{0\}\sub \K^{d_n}$).
Define $P_n:=\Delta^{d_n}_1$,
$Q_n:=\psi_n(P_n)$, and $\phi_n:=\psi_n|_{P_n}^{Q_n}\!:
P_n\to Q_n$ for $n\geq \ell$.
Then $P_\phi:=\bigcup_{n\in \N} P_n$
is open in~$\K^s$, $Q_\phi:=\bigcup_{n\in \N}Q_n$
is open in~$M$,
and $\phi:=\dl\, \phi_n\!: P_\phi\to Q_\phi$\vspace{-.8mm}
is an element of $\cA$, with $x\in Q_\phi$,
as desired.\vspace{1mm}

{\em Compatibility of the charts.}
Assume that $\phi:=\dl\,\phi_n\!: P_\phi\to Q_\phi$\vspace{-.8mm}
and $\psi:=\dl\, \psi_n\!: P_\psi\to Q_\psi$
are elements of $\cA$,
where $\phi_n\!: P_n\to Q_n$
and $\psi_n\!: A_n\to B_n$.
Suppose that $x\in \phi^{-1}(Q_\psi)$.
Then $\phi(x)\in Q_\phi\cap Q_\psi$,
entailing that there exists
$\ell\in \N$ such that $\phi(x)\in Q_\ell\cap B_\ell$.
Then $x\in P_n\cap\phi_n^{-1}(B_n)=:X_n$
for all $n\geq \ell$.
Since $X_n$ is open in $\K^{d_n}$
and $X_\ell\sub X_{\ell+1}\sub \cdots$,
the union $X:=\bigcup_{n\geq \ell}X_n$
is open in~$\K^s$.
Furthermore, the coordinate changes
$\tau_n:=\psi_n^{-1}|_{Q_n\cap B_n}\circ \phi_n|_{X_n}\!:
X_n\to \psi_n^{-1}(Q_n)=:Y_n$
are $C^\infty_\K$-diffeomorphisms,
for all $n\geq \ell$.
By Lemma~\ref{silva}\,(a),
the map $\psi^{-1}|_{\phi(X)}^Y\circ \phi|_X
=\dl_{n\geq \ell}\, \tau_n\!: X\to \bigcup_{n\geq \ell}Y_n=:Y$
is $C^\infty_\K$,
entailing that the bijection $\tau:=
\psi^{-1}|_{Q_\phi\cap Q_\psi}\circ \phi|_{\phi^{-1}(Q_\psi)}\!:
\phi^{-1}(Q_\psi)\to \psi^{-1}(Q_\phi)$
is $C^\infty_\K$
on some open neighbourhood of~$x$.
As $x$ was arbitrary, $\tau$ is $C^\infty_\K$
and the same reasoning shows that so is $\tau^{-1}$.
Thus $\cA$ is an atlas making $M$ a $C^\infty_\K$-manifold
modelled on~$\K^s$.\vspace{1mm}

{\em Each $\lambda_n$ is smooth.}
To see this, assume that $n\in \N$ and $x\in M_n$.
As just shown,
there exists a chart $\phi\!: P_\phi\to Q_\phi$
in $\cA$,
say $\phi=\dl\, \phi_k$\vspace{-.8mm}
with charts $\phi_k\!: P_k\to Q_k\sub M_k$
for $k\in \N$,
such that $x\in P_n$.
Then $\phi^{-1}\circ \lambda_n\circ
\phi_n=\phi^{-1}\circ \phi_n\!: \K^{d_n}\supseteq P_n\to P\sub \K^s$
is the inclusion map and hence smooth,
and its differential at $x$ is injective.
Hence $\lambda_n$ is smooth on the open neighbourhood
$Q_n$ of~$x$,
and $T_x(\lambda_n)$
is injective.
As $x$ was arbitrary, $\lambda_n$ is smooth.\vspace{1mm}

{\em Direct limit property and uniqueness.}
Fix $r\in \N_0\cup\{\infty\}$.
Assume that $Y$ is a $C^r_\K$-manifold
modelled on a topological $\K$-vector space~$E$
and $f_n\!: M_n\to Y$
a $C^r_\K$-map for each $n\in \N$ such that $(Y,(f_n)_{n\in \N})$
is a cone over $\sys$;
thus $f_n|_{M_m}=f_m$ if $n\geq m$.
Then there is a uniquely determined
map $f\!: M\to Y$ such that $f|_{M_n}=f_n$
for all $n\in \N$.
Since $M=\dl\,M_n$\vspace{-.8mm}
as a topological space,
$f$ is continuous.
If $x\in M$,
we find a chart $\phi\!: P_\phi\to Q_\phi$ of $M$
around~$x$
in the atlas $\cA$, where $\phi=\dl \,\phi_n$\vspace{-.8mm}
for charts $\phi_n\!: P_n\to Q_n\sub M_n$.
Let $\psi\!: V\to W\sub Y$ be a chart for~$Y$,
where $V\sub E$ is open.
Then $U:=(f\circ \phi)^{-1}(W)$
is an open subset of $P_\phi\sub \K^s$,
and $U_n:=U\cap P_n$ is open in $P_n\sub \K^{d_n}$
for each~$n$.
Consider
$g:=\psi^{-1}\circ (f\circ \phi)|_U^W\!: U\to V$.
Then $g|_{U_n}=\psi^{-1}\circ (f_n\circ \phi_n)|_{U_n}^W\!:
U_n\to V$ is $C^r_\K$ for each $n\in \N$.
Hence $g$ is $C^r_\K$
by Lemma~\ref{silva}\,(a), whence so is $f$
on the open neighbourhood $Q_\phi$ of $x$
and hence on all of~$M$, as $x$ was arbitrary.
Thus
$(M,(\lambda_n)_{n\in \N})=\dl\,\sys$\vspace{-.8mm}
in the category of $C^r_\K$-manifolds,
for all $r\in \N_0\cup\{\infty\}$.
The uniqueness of a $C^\infty_\K$-manifold structure
on~$M$ with the described properties
follows from the universal
property of direct limits.
\end{proof}
{\bf Convention.}
Throughout the remainder of this section,
$M$ will be
equipped with the $C^\infty_\K$-manifold structure
just defined.
In the proofs,
we shall always
reduce to the case where $I=\N$ and $M_1\sub M_2\sub \ldots$
(by the above argument), without further mention.
\begin{prop}\label{mfdh}
If $\,\F\sub \K$ is a non-discrete, closed
subfield, then $(M,(\lambda_i)_{i\in I})=\dl\, \sys$\vspace{-.8mm}
also in the category of $C^\infty_\F$-manifolds.
$($E.g., $\K=\C$, $\F=\R\,)$.
\end{prop}
\begin{proof}
Let $\cA$ be the $C^\infty_\K$-atlas of $M$ described
in the proof of Theorem~\ref{createmfd}.
Given a non-discrete closed
subfield $\F\sub \K$, let 
$\cA_\F$ be the corresponding atlas
obtained when considering each $M_i$
merely as a $C^\infty_\F$-manifold over $\F$.
Then $\cA\sub \cA_\F$,
entailing that $(M,(\lambda_n)_{n\in N})=\dl\, \sys$\vspace{-1mm}
also in the category of $C^\infty_\F$-manifolds.
\end{proof}
\begin{prop}\label{mfdb}
Assume that $U_i\sub M_i$
is open and $\lambda_{i,j}(U_j)\sub U_i$
whenever $i\geq j$.
Then $U:=\bigcup_{i\in I} U_i$ is open in~$M$.
For the $C^\infty_\K$-manifold
structure induced by $M$ on its open subset~$U$,
we have
$(U,(\lambda_i|_{U_i}^U)_{i\in I})=\dl
((U_i)_{i\in I}, (\lambda_{i,j}|_{U_j}^{U_i})_{i\geq j})$\vspace{-1mm}
in the category of $C^\infty_\K$-manifolds.
\end{prop}
\begin{proof}
Given open subsets $U_n\sub M_n$
such that $M_1\sub M_2\sub \cdots$,
their union $U:=\bigcup_{n\in \N}U_n$
is open in $M$
and the topology induced by $M$ on~$U$
makes $U$ the direct limit $\dl\, U_n$\vspace{-.8mm}
(Lemma~\ref{base}\,(b)). We define an atlas
$\cA_U$ for $U$ turning $U$
into the direct limit of the $C^\infty_\K$-manifolds
$U_n$, analogous to the definition of $\cA$ in the proof
of~(a). Then $\cA_U\sub \cA$, whence
$(U,\cA_U)$ coincides with $U$, considered as an open
submanifold of~$M$.
\end{proof}
\begin{prop}\label{mfdc}
Assume that
$f\!: X \to M$ is a $C^r_\K$-map,
where $r\in \N_0\cup\{\infty\}$
and $X$ is a
$C^r_\K$-manifold modelled
on a metrizable topological $\K$-vector space~$E$
$($\/or a metrizable, locally path-connected topological space,
if $r=0$\/$)$.
Then every $x\in X$ has an open neighbourhood
$S$ such that $f(S)\sub \lambda_i(M_i)$
for some $i\in \N$ and such that
$\lambda_i^{-1}\circ f|_S^{\lambda_i(M_i)}\!: S\to M_i$
is $C^r_\K$.
\end{prop}
\begin{proof}
Let $x\in X$. The assertion
being local, in the case of manifolds
we may assume that $X$ is an open subset of~$E$.
Choose a metric
$d$ on $X$ defining its topology,
and $k\in \N$ such that $f(x)\in M_k$.
Let $\phi=\dl\, \phi_n\!: P\to Q$\vspace{-1mm}
be a chart of~$M$ around $f(x)$,
where $\phi_n$
is a chart of $M_n$ for all $n\geq k$,
of the form $\phi_n\!: \Delta^{d_n}_1\to Q_n\sub M_n$
(see proof of Theorem~\ref{createmfd}).
If $f^{-1}(Q_n)$ is not a neighbourhood
of~$x$ for any $n\geq k$,
we find $x_n\in f^{-1}(Q) \setminus f^{-1}(Q_n)$
such that $d(x_n,x)<2^{-n}$.
Thus $x_n\to x$, entailing that
$C:=\{f(x_n)\!: n\in \N\}\cup\{f(x)\}$
is a compact subset of $Q$
such that $C\not\sub Q_n$
for any $n\geq k$. Since $Q=\dl\, Q_n$,\vspace{-.8mm}
this contradicts Lemma~\ref{base}\,(d).
Hence there exists $n\geq k$ such that $f^{-1}(Q_n)$
is a neighbourhood of~$x$. Let $S:=(f^{-1}(Q_n))^0$
be its interior. Then $S\to \K^s$,
$y\mto \phi^{-1}(f(y))=\phi_n^{-1}(f(y))$
is a $C^r_\K$-map taking its values in
the closed vector subspace $\K^{d_n}$ of $\K^s$,
whence also its co-restriction $\phi_n^{-1}\circ
f|_S^{Q_n}\!: S\to \Delta^{d_n}_1$ is $C^r_\K$
\cite[La.\,10.1]{BGN}.
As $\phi_n$ is a chart,
this means that
$f|_S^{M_n}$ is $C^r_\K$.
\end{proof}
\begin{prop}\label{mfdd}
If $\K\in \{\R,\C\}$ and $x\in M$, where $x=\lambda_j(y)$ say,
then the connected component $C$ of $x\in M$
in~$M$ is $\bigcup_{i\in I}\lambda_i(C_i)\isom \dl_{i\geq j}
\,C_i$,\vspace{-2mm}
where $C_i$ is the connected component of $\lambda_{i,j}(y)$ in $M_i$.
\end{prop}
\begin{proof}
Given $x\in M_n$,
we let $C$ and $C_m$ be its connected component in~$M$
and $M_m$, respectively,
for $m\geq n$.
Then
$\bigcup_{m\geq n}C_m\sub C$.
If $z\in C$, then we find a continuous curve
$\gamma\!: [0,1]\to M$ such that
$\gamma(0)=x$ and $\gamma(1)=z$.
Since $[0,1]$ is compact,
using Proposition~\ref{mfdc} we find $m\geq n$
such that $\gamma([0,1])\sub M_m$,
and such that $\gamma|^{M_m}\!: [0,1]\to M_m$
is continuous. Thus $z\in C_m$.
Hence indeed $C=\bigcup_{m\geq n}C_m$.
\end{proof}
\begin{prop}\label{mfde}
If $\K$ is $\R$ or $\C$
and $M_i$ is paracompact for
each $i\in I$, then
the $C^\infty_\K$-manifold $M$
also is a
$c^\infty_\K$-manifold,
and $(M,(\lambda_i)_{i\in I})\!=\!
\dl\, \sys$\vspace{-.8mm} in the category of
$c^\infty_\K$-manifolds.
Furthermore,
$M$ is smoothly paracompact
as a $C^\infty_\R$-manifold{\rm :}
For
every open cover of $M$,
there exists a $C^\infty_\R$-partition
of unity subordinate to the cover.
\end{prop}
\begin{proof}
Assume that $\K=\R$.
In order that $M$ be smoothly paracompact,
we only need to show that every connected component
$C$ of~$M$
is smoothly paracompact.
Pick $c\in C$.
We may assume that $c\in M_1$
after passing to a cofinal subsystem;
we
let $C_n$ the connected component
of $c$ in $M_n$
for each~$n\in \N$.
Then $\bigcup_{n\in \N}C_n$
is the connected component of
$c$ in~$M$
(see Proposition~\ref{mfdd}) and
hence coincides with~$C$;
furthermore, $C=\dl\, C_n$,\vspace{-.8mm} by
Proposition~\ref{mfdb}.
After replacing $M$ with $C$
and $M_n$ with $C_n$
for each~$n$, we may assume that each $M_n$
is a connected, paracompact finite-dimensional
$C^\infty_\R$-manifold
and hence $\sigma$-compact.
This entails that $M=\bigcup_{n\in \N} M_n$
is $\sigma$-compact and therefore Lindel\"{o}f.
Hence, by \cite[Thm.\,16.10]{KaM},
$M$ will be smoothly paracompact
if we can show that $M$ is {\em smoothly
regular\/} in the sense
that, for every
$x\in M$ and open neighbourhood
$\Omega$ of $x$ in~$M$,
there exists a smooth function
(``bump function'')
$f\!: M\to \R$ such that $f(x)\not=0$
and $f|_{M\setminus \Omega}=0$.\vspace{1mm}

If each $\lambda_{n,m}$ is a topological
embedding onto a closed
submanifold, then $M$ is a
regular topological space
(see \cite[Prop.\,4.3\,(ii)]{Han}),
whence smooth regularity
passes from the modelling space\footnote{
See \cite[Thm.\,16.10]{KaM}
or \cite{DIR}, proof of Thm.\,6.4
for the smooth regularity of $\R^\infty$.}
to~$M$ (cf.\ \cite[proof of Thm.\,6.4]{DIR}).
In the fully general case to be investigated here,
we do not know {\em a priori\/}
that $M$ is regular,
whence we have to prove smooth regularity
of~$M$ by hand.
Essentially,
we need to go once more through our construction of charts
and build up bump functions step by step.
Let $x\in M$ and $\Omega$ be an open neighbourhood
of~$x$ in~$M$. Passing to a co-final subsequence,
we may assume that $x\in M_1$.\vspace{1mm}

Let $r_n:=1+2^{-n}$ for $n\in \N$
and $\Delta^{d_n}_{r_n}:=\Delta^{d_n}_{r_n}(\R)$.
Pick a chart
$\psi_1\!: \Delta^{d_1}_{r_1}\to W_1\sub
M_1$ of $M_1$ around~$x$, such that $\psi_1(0)=x$
and such that $W_1$ is relatively compact in~$M_1\cap \Omega$.
Define $Q_1:=\psi_1(\Delta^{d_1}_1)$.
We choose compact subsets $K_{1,j}$ of~$M_1$ such that
$W_1\sub K_{1,1}^0\sub K_{1,1}\sub K_{1,2}^0\sub K_{1,2}\sub K_{1,3}^0\sub
\cdots$
and $M_1=\bigcup_{j\in \N}K_{1,j}$.
There exists a smooth function
$h_1\!: \Delta^{d_1}_{r_1}\to \R$ such that
$\Supp(h_1)\sub \Delta^{d_1}_1$
and $h_1(0)=1$.
Define $f_1\!: M_1\to \R$,
$f_1(y):=0$ if $y\not\in W_1$,
$f_1(y):=h_1(\psi_1^{-1}(y))$ if $y\in W_1$.
Then $f_1$ is smooth, $\Supp(f_1)\sub Q_1$,
and $f_1(x)=1$.\vspace{1mm}

The Extension Lemma~\ref{extpoly}
provides a chart $\psi_2\!: \Delta^{d_2}_{r_2} \to W_2\sub M_2$
onto an open, relatively compact subset
$W_2$ of $M_2\cap \Omega$
such that $\psi_2|_{\Delta^{d_1}_{r_2}}=
\psi_1|_{\Delta^{d_1}_{r_2}}$.
We choose compact subsets $K_{2,j}$ of~$M_2$ such that
$K_{1,j}\sub K_{2,j}^0$
and
$W_2\sub K_{2,1}^0\sub K_{2,1}\sub K_{2,2}^0\sub K_{2,2}\sub
\cdots$
and $M_2=\bigcup_{j\in \N}K_{2,j}$.
Then $K_{1,1}\setminus Q_1$
is a compact subset of $M_1$ and hence also
of $M_2$. Therefore
$A:=\psi_2^{-1}(K_{1,1}\setminus Q_1)$ is closed in
$\Delta^{d_2}_{r_2}$,
and it does not meet the compact subset
$\Supp(h_1)\sub \Delta^{d_1}_1\sub \Delta^{d_2}_{r_2}$
(which is mapped into $Q_1$ by $\psi_2$).
Hence, there exists $\ve\in \;]0,1[$
such that $A\cap (\Supp(h_1)\times \Delta^{d_2-d_1}_\ve)=\emptyset$.
We let $\xi\!: \R\to \R$
be a smooth function such that
$\xi(0)=1$
and $\Supp(\xi)\sub \;]{-\ve^2},\ve^2[$.
Then
\[
h_2\!: \Delta^{d_2}_{r_2}\to\R,\qquad
h_2(y,z):=h_1(y)\cdot \xi((\|z\|_2)^2)
\quad \mbox{for $y\in \Delta^{d_1}_{r_2}$,
$z\in \Delta^{d_2-d_1}_{r_2}$}
\]
(where $\|.\|_2$ is the euclidean norm
on $\R^{d_2-d_1}$)
is a smooth map such that $\Supp(h_2)\sub
\Delta^{d_2}_1$.
Then $f_2(y):=0$ if $y\not\in W_2$,
$f_2(y):=h_2(\psi_2^{-1}(y))$ for $y\in W_2$
defines a smooth function
$f_2\!: M_2\to \R$.
We have $\Supp(f_2)\sub Q_2:=\psi_2(\Delta^{d_2}_1)$,
and $f_2|_{K_{1,1}}=f_1|_{K_{1,1}}$,
because $f_2|_{Q_1}=f_1|_{Q_1}$
by definition
and also
$f_2|_{K_{1,1}\setminus Q_1}=0=f_1|_{K_{1,1}\setminus Q_1}$.\vspace{1mm}

Proceeding in this way, we
find charts $\psi_n\!: \Delta_{r_n}^{d_n}\to W_n\sub M_n$
with relatively compact image
$W_n\sub \Omega$,
compact subsets $K_{n,j}$
of $M_n$ with union~$M_n$ such that
$W_n\sub K_{n,1}$, $K_{n,j}\sub K_{n,j+1}^0$
and $K_{n-1,j}\sub K_{n,j}^0$ for all $n,j\in \N$,
$n\geq 2$;
and smooth maps $f_n\!: M_n\to \R$
such that $\Supp(f_n)\sub Q_n:=\psi_n(\Delta^{d_n}_1)$
and
$f_{n+1}|_{K_{n,n}}=f_n|_{K_{n,n}}$
for all $n\in \N$, whence
\begin{equation}\label{maintrick}
f_m|_{K_{n,n}}=f_n|_{K_{n,n}}
\qquad \mbox{for all $n,m\in \N$ such that $m\geq n$.}
\end{equation}
Let $U_n$ be the interior $K_{n,n}^0$
of $K_{n,n}$ in $M_n$. Then
$U_1\sub U_2\sub \cdots$ and
$M=\bigcup_{n\in \N}U_n$, whence
$M=\dl\, U_n$\vspace{-.8mm} as a smooth manifold
by Proposition~\ref{mfdb}. By (\ref{maintrick}),
the smooth maps $f_n|_{U_n}$
form a cone and hence
induce a smooth map $f\!: M\to \R$,
such that $f|_{U_n}=f_n|_{U_n}$ for each $n\in \N$.
Then $f(x)=f_1(x)=1$.
If $y\in M\setminus \Omega$,
we find $n\in \N$ such that $y\in U_n$.
Then $f(y)=f_n(y)=0$ because
$\Supp(f_n)\sub Q_n\sub W_n\sub \Omega$.
Hence $f$ is a bump function
around~$x$ carried by~$\Omega$,
as desired. Thus $M$ is smoothly paracompact.\vspace{1mm}

{\em Direct limit properties when $\K\in \{\R,\C\}$.}
Since $M$ (resp., its underlying real
manifold) is smoothly regular,
$M$ is {\em smoothly Hausdorff},
i.e., $C^\infty_\R(M,\R)$
separates points on~$M$. Let $\cA$ be an $C^\infty_\K$-atlas
for~$M$.
Being a smoothly Hausdorff $C^\infty_\K$-manifold
modelled on a Mackey complete locally convex
space, $(M,\cA)$ can be made a $c^\infty_\K$-manifold
$(c^\infty(M),\cA)$ by replacing its topology
with the final
topology with respect to
the given charts, when the
topology on the modelling space
has been replaced with its $c^\infty$-topology.
Since $c^\infty(\K^s)=\K^s$,
the topology on~$M$ remains unchanged,
and thus $c^\infty(M)=M$. In view of
Lemma~\ref{silva}\,(b),
the desired direct limit properties can be checked as in the proof
of Theorem~\ref{createmfd}.
\end{proof}
\begin{prop}\label{mfdf}
Assume that also
$\cT:=((N_i)_{i\in I}, (\mu_{i,j})_{i\geq j})$
is a direct system of finite-dimensional $C^\infty_\K$-manifolds
and injective $C^\infty_\K$-immersions,
over the same index set.
Then also $\cP:=
((M_i\times N_i)_{i\in I}, (\lambda_{i,j}\times\mu_{i,j})_{i\geq j})$
is such a direct system.
Let $(N,(\mu_i))=\dl\,\cT$.\vspace{-.8mm}
The $C^\infty_\K$-maps $\eta_i:=\lambda_i\times \mu_i\!:
M_i\times N_i\to M\times N$
define a cone $(M\times N,(\eta_i)_{i\in I})$ over $\cP$,
which
induces a $C^\infty_\K$-diffeomorphism
$\eta \!: \dl \,(M_i\times N_i)\to \bigl(\dl \,M_i\bigr)\times
\bigl(\dl\, N_i\bigr)$.\vspace{-1mm}
\end{prop}
\begin{proof}
Let $e_n:=\dim_\K(N_n)$ and $t:=\sup\{e_n\!: n\in \N\}$.
The natural map $\zeta\!: \K^{s+t}=\dl\, \K^{d_n+e_n}\to \K^s\times
\K^t$\vspace{-.8mm} analogous to $\eta$ is an isomorphism
of topological vector spaces
(\cite[Prop.\,3.3]{DIR}; \cite[Thm.\,4.1]{HSTH}).
Let $\cA$ be the atlas for $M=\bigcup_{n\in \N}M_n$
from the proof of Theorem~\ref{createmfd}; let
$\cB$ and $\cC$ be analogous
atlases for $N=\bigcup_{n\in \N}N_n$
and $P:=\bigcup_{n\in \N}(M_n\!\times \!N_n)$.
Then $\cD:=\{\phi\times \psi\!: \phi\in \cA,\psi\in \cB\}$
is a $C^\infty_\K$-atlas making
$M\times N$ the direct product of~$M$ and $N$
in the category of $C^\infty_\K$-manifolds.
Since $\{(\phi\times \psi)\circ \zeta|_{\zeta^{-1}(P_\phi\times P_\psi)}\!:
\phi\in \cA,\psi\in \cB\}\sub
\cC$, the map
$\eta=\id\!: P\to M\times N$
is a $C^\infty_\K$-diffeomorphism.
\end{proof}
\begin{prop}\label{mfdg}
If $\K=\R$, each $M_i$ is a finite-dimensional,
real analytic
manifold and each $\lambda_{i,j}$ an injective,
real analytic immersion, then
there exists a $c^\omega_\R$-manifold structure on $M$ such that
$M=\dl\, \sys$\vspace{-.8mm}
in the category of
$c^\omega_\R$-manifolds $($and $c^\omega_\R$-maps$)$,
and which is compatible with
the above $C^\infty_\R$-manifold structure on~$M$.
Analogues of Propositions
{\/\rm \ref{mfdb}, \ref{mfdc}\/} and {\/\rm \ref{mfdf}\/}
hold for the $c^\omega_\R$-manifold structures
$($and $c^\omega_\R$-maps$)$.
\end{prop}
\begin{proof}
Using the $C^\omega_\R$-case of the Extension Lemma~\ref{extpoly},
the construction described in the proof of
Theorem~\ref{createmfd}
provides a subatlas $\cB$ of the $C^\infty_\R$-atlas
$\cA$ consisting of charts $\phi=\dl\,\phi_n$\vspace{-.8mm}
where each $\phi_n$ is a
$C^\omega_\R$-diffeomorphism.
Using Lemma~\ref{silva}\,(b), the above arguments show
that the chart changes
for charts in $\cB$ are~$c^\omega_\R$,
whence indeed $M$ has a compatible
$c^\omega_\R$-manifold structure.
Similarly, Lemma~\ref{silva}\,(b)
entails the validity of the
$c^\omega_\R$-analogues of Proposition~\ref{mfdb} and \ref{mfdf}. The proof of
Proposition~\ref{mfdc}
carries over directly.
\end{proof}
\begin{rem}\label{homotopy}
For the singular homology groups of~$M$
over a commutative ring~$R$,
we have
$(H_k(M),\,(H_k(\lambda_i))_{i\in I})=
\dl\, (H_k(M_i),\, (H_k(\lambda_{i,j})_{i\geq j}))$
for all $k\in \N_0$,
as a consequence of Proposition~\ref{mfdc}
(and Proposition~\ref{mfdd}).
Likewise, given $\ell\in I$
and $x\in M_\ell$,
for the homotopy groups we have
$\pi_k(M, \lambda_\ell(x))=
\dl\, \pi_k(M_i,\lambda_{i,\ell}(x))$.
\end{rem}
\section{Direct limits of finite-dimensional Lie groups}\label{secgps}
Let $\K$ be $\R$, $\C$ or a local field,
and $\sys:=((G_i)_{i\in I},(\lambda_{i,j})_{i\geq j})$
be a countable direct system of finite-dimensional
$C^\infty_\K$-Lie groups~$G_i$
and $C^\infty_\K$-homomorphisms $\lambda_{i,j}\!: G_j\to G_i$;
if $\car(\K)>0$, we assume that each $\lambda_{i,j}$
is an injective immersion.
In this section,
we construct a direct limit Lie group for~$\sys$,
and discuss some of its properties.
\begin{rem}
(a)
As in the classical real
and complex cases, also
every $C^\infty_\K$-Lie group over a local
field~$\K$ of characteristic~$0$
admits a $C^\infty_\K$-compatible analytic
Lie group structure, and every $C^\infty_\K$-homomorphism
between such groups is $\K$-analytic~\cite{ANA}.\vspace{1mm}

(b) Note that
the squaring map $\sigma\!: \F_2[[X]]^\times\to \F_2[[X]]^\times$,
$\sigma(x):=x^2$ is an analytic (and hence smooth)
homomorphism which is injective
(since $\F[[X]]^\times$
is isomorphic to the power
$(\Z_2)^\N$ of the group of $2$-adic integers
as a group by \cite[Chap.\,II-4, Prop.\,10]{Wei},
and thus torsion-free)
but {\em not\/}
an immersion because $f'(x)=2\, \id=0$ for all $x\in \F_2[[X]]^\times$
(where $\F_2[[X]]$ denotes the field of formal Laurent series
over the field $\F_2:=\Z/2\Z$).
This explains that an extra hypothesis is needed
in positive characteristic.\footnote{Slightly more generally,
to establish Theorem~\ref{limgps} for $\car(\K)>0$,
it would be enough to assume that
that $G_j/N_j$ admits a $C^\infty_\K$-Lie group structure
for each $j\in I$
which makes the quotient map $G_j\to G_j/N_j$ a submersion,
where $N_j:=\wb{\bigcup_{i\geq j}\ker\, \lambda_{i,j}}$,
and that the induced homomorphisms
$G_j/N_j\to G_i/N_i$
be immersions, for all $i\geq j$.}
\end{rem}
{\bf Associated injective quotient system.}
If $\K =\R$
(or if $\car(\K)>0$, in which case we obtain trivial
subgroups),
we let
$N_j:=\wb{\bigcup_{i\geq j}\ker\, \lambda_{i,j}}$
for $j\in I$.
If
$\K$ is a local field of characteristic~$0$,
we let $N_j=
\bigcup_{i\geq j}\ker\, \lambda_{i,j}$
for $j\in I$
and note that $N_j$
is locally closed and hence closed in $G_j$,
because $G_j$ has an open compact subgroup~$U$
which satisfies an ascending chain condition
on closed subgroups.
If $\K=\C$, we let $N_j\sub G_j$
be the
intersection of all closed
complex Lie subgroups~$S$
of $G_j$ such that $\bigcup_{i\geq j}\ker\, \lambda_{i,j}\sub S$.
Then $N_j$ is a complex Lie subgroup
of $G_j$ (as $G_j$ satisfies a descending
chain condition on closed, connected subgroups),
and thus $N_j$ is the smallest
closed, complex Lie subgroup of~$G_j$
containing
$\bigcup_{i\geq j}\ker\, \lambda_{i,j}$.
By minimality, $N_j$
is invariant under inner automorphisms and hence
a normal subgroup of~$G_j$.

Then, in either case,
there is a uniquely
determined $\K$-Lie group structure on $\wb{G}_j:=G_j/N_j$
which makes the canonical quotient homomorphism
$q_j\!: G_j\to \wb{G}_j$ a submersion.
Each $\lambda_{i,j}$ factors to a $C^\infty_\K$-homomorphism
$\wb{\lambda}_{i,j}\!: \wb{G}_j\to \wb{G}_i$,
uniquely
determined by $\wb{\lambda}_{i,j}\circ q_j=q_i\circ \lambda_{i,j}$.
Then
$\wb{\sys}=((\wb{G}_i)_{i\in I}, (\wb{\lambda}_{i,j})_{i\geq j})$
is a direct system of finite-dimensional $C^\infty_\K$-Lie groups
and {\em injective\/} $C^\infty_\K$-homomorphisms
$\wb{\lambda}_{i,j}\!: \wb{G}_j\to \wb{G}_i$;
it is called the {\em injective quotient system
associated with $\sys$\/}
(cf.\ \cite{NRW1}). Each $\wb{\lambda}_{i,j}$
is an injective immersion of class $C^\infty_\K$.
\begin{rem}\label{preparecx}
The example
$\C^\times\stackrel{\sigma}{\to}\C^\times \stackrel{\sigma}{\to}\cdots$
with the squaring map $\sigma\!: \C^\times\to\C^\times$,
$\sigma(z):=z^2$ shows that $\wb{\bigcup_{i\geq j}\ker\, \lambda_{i,j}}$
need not be a {\em complex\/} Lie subgroup
of $G_j$.
\end{rem}
$\;$\vspace{-4mm}\pagebreak

\begin{thm}\label{limgps}
For $\sys$ and $\wb{\sys}$ as before,
the following holds:\vspace{-1mm}
\begin{itemize}
\item[{\rm (a)}]
A direct limit $(G,(\wb{\lambda}_i)_{i\in I})$
for $\wb{\sys}$
exists in the category of $C^\infty_\K$-Lie groups
modelled on topological $\K$-vector spaces;
$G$ is modelled on the locally convex topological
$\K$-vector space $\K^s$, where
$s:=\sup\{\dim_\K \wb{G}_i\!: i\in I\}\in \N_0\cup\{\infty\}$.
Forgetting the $\K$-Lie group structure,
we have $(G,(\wb{\lambda}_i)_{i\in I})=\dl\, \wb{\sys}$\vspace{-.8mm}
also in the categories
of sets, abstract groups,
topological spaces,
topological groups,
the category of
$C^\infty_\K$-manifolds modelled on topological
$\K$-vector spaces,
and the category of $C^\infty_\F$-Lie groups
modelled on topological $\F$-vector spaces,
for every non-discrete closed subfield~$\,\F$ of~$\,\K$.
Furthermore, $L(\wb{\lambda}_i)\!: L(\wb{G}_i)\to L(G)$
is injective for each $i\in I$, and
\begin{equation}\label{descLie}
\bigl(L(G),(L(\wb{\lambda}_i))_{i\in I}\bigr)\;
=\; \dl\, \bigl((L(\wb{G}_i))_{i\in I},
(L(\wb{\lambda}_{i,j}))_{i\geq j}\bigr)\vspace{-1mm}
\end{equation}
in the category of locally convex
$\K$-vector spaces
$($\/and in the categories
of sets, $\K$-Lie algebras, topological spaces,
topological $\K$-Lie algebras,
topological $\K$-vector spaces,
and $C^\infty_\K$-manifolds; also
in the category of $c^\omega_\K$-manifolds,
if $\K\in\{\R,\C\}$\/$)$.
\item[\rm (b)]
Set $\lambda_i:=\wb{\lambda}_i\circ q_i$
for $i\in I$. If $\K\not=\C$
or if $\,\K=\C$ and $N_j=\wb{\bigcup_{i\geq j}\ker\lambda_{i,j}}$
for each $j\in I$,
then
$(G,(\lambda_i)_{i\in I})=\dl\, \sys$\vspace{-.8mm}
in the category of $C^\infty_\K$-Lie groups
modelled on topological $\K$-vector spaces,
and also in the categories of
smooth manifolds modelled on topological
$\K$-vector spaces,
Hausdorff topological spaces,
and $($Hausdorff\,$)$ topological groups.
\item[\rm (c)]
If $\K=\C$,
then
$(G,(\lambda_i)_{i\in I})=\dl\, \sys$\vspace{-.8mm}
in the category of complex Lie groups
modelled on complex locally convex spaces.
\item[\rm (d)]
If $\,\K=\R$,
then $G$ is
a $c^\infty_\R$-regular $c^\infty_\R$-Lie group
which is smoothly paracompact.
Furthermore, $(G,(\lambda_i)_{i\in I})=\dl\, \sys$\vspace{-1mm} in the category
of $c^\infty_\R$-Lie groups,
and in the category of $c^\infty_\R$-manifolds.
\item[\rm (e)]
If $\,\K=\C$, then $G$ is a $c^\infty_\C$-regular,
$c^\infty_\C$-Lie group such that
$(G,(\lambda_i)_{i\in I})=\dl\,\sys$\vspace{-1mm}
in the category of $c^\infty_\C$-Lie
groups and
$(G,(\lambda_i)_{i\in I})=\dl\, \wb{\sys}$\vspace{-1mm}
in the category of $c^\infty_\C$-manifolds.
\item[\rm (f)]
If $\,\K=\R$, then there exists a $c^\omega_\R$-regular,
$c^\omega_\R$-Lie group structure on $G$,
compatible with the $C^\infty_\R$-Lie group structure from
{\rm (a)},
such that
$(G,(\lambda_i)_{i\in I})=\dl\, \sys$\vspace{-.8mm}
in the category of $c^\omega_\R$-Lie groups.
For the underlying
$c^\omega_\R$-
manifold, we have $(G,(\lambda_i)_{i\in I})=\dl\, \sys$\vspace{-1mm}
in the category of such manifolds.
\end{itemize}
\end{thm}
\begin{proof}
(a) Let $(G,(\wb{\lambda}_i)_{i\in I})$
be a direct limit for $\wb{\sys}$
in the category of topological groups;
then $(G,(\wb{\lambda}_i)_{i\in I})=\dl\, \wb{\sys}$\vspace{-1mm}
also in the categories of sets, groups,
and topological spaces \cite[Thm.\,2.7]{TSH}.
Thus Theorem~\ref{createmfd}
provides a $C^\infty_\K$-manifold structure
on~$G$ making $(G,(\wb{\lambda}_i)_{i\in I})$
a direct limit of $\wb{\sys}$
in the category of $C^\infty_\K$-manifolds modelled
on topological $\K$-vector spaces,
and also in the category of $C^\infty_\F$-manifolds,
for every non-discrete closed subfield $\F\sub \K$
(Proposition~\ref{mfdh}).
Let $\theta \!: G\to G$, $g\mto g^{-1}$ and $\theta_i\!: \wb{G}_i\to \wb{G}_i$
be the inversion maps;
let $\mu\!: G\times G$ and $\mu_i\!: \wb{G}_i\times \wb{G}_i\to \wb{G}_i$.
Then $\theta:=\dl\, \theta_i$\vspace{-.8mm}
is $C^\infty_\K$,
as $G=\dl\, \wb{G}_i$\vspace{-.8mm}
in the category of $C^\infty_\K$-manifolds.
Likewise,
$\mu=\dl\,\mu_i$\vspace{-.8mm}
is $C^\infty_\K$
on $\dl\, (\wb{G}_i\times \wb{G}_i)$\vspace{-.8mm}
and hence on $G\times G$, in view of
Proposition~\ref{mfdf}.
Hence $G$ is a $C^\infty_\K$-Lie group.
Every cone $(H,(f_i)_{i\in I})$
of $C^\infty_\K$-homomorphisms $f_i\!: \wb{G}_i\to H$
to a $C^\infty_\K$-Lie group~$H$
uniquely determines a map
$f\!: G\to H$ such that $f\circ \wb{\lambda}_i=f_i$
for all~$i$;
then $f$ is a homomorphism since $G=\dl\, \wb{G}_i$\vspace{-.8mm}
as a group, and $f$ is $C^\infty_\K$
since $G=\dl\, \wb{G}_i$\vspace{-.8mm}
as a $C^\infty_\K$-manifold.
Thus $G=\dl\, \wb{G}_i$
as a $C^\infty_\K$-Lie group
(and, likewise, as a $C^\infty_\F$-Lie group).\vspace{1mm}

{\em Determination of the Lie algebra of $L(G)$.}
Since $L(\wb{\lambda}_{i,j})$
is injective for all $i\geq j$,
$\cT:=((L(\wb{G}_i))_{i\in I}, (L(\wb{\lambda}_{i,j}))_{i\geq j})$
is a countable, strict direct system of
Lie algebras. We recall
from \cite{DIR} or \cite{HSTH} that $\cT$
has a direct limit $(\g,(\phi_i)_{i\in I})$
in the category of topological $\K$-Lie algebras;
here $\g$ carries the finite
topology (see {\bf \ref{new1.8}}), each $\phi_i$ is injective,
and
$(\g,(\phi_i)_{i\in I})=\dl\, \cT$\vspace{-.8mm}
also holds in the categories
of sets, $\K$-Lie algebras, topological spaces,
topological $\K$-vector spaces,
and locally convex topological $\K$-vector spaces.
By Lemma~\ref{silva}\,(a),
furthermore
$(\g,(\phi_i)_{i\in I})=\dl\, \cT$\vspace{-.8mm}
in the category of $C^\infty_\K$-manifolds
and $C^\infty_\K$-maps
(and also in the category of
$c^\omega_\K$-manifolds and $c^\omega_\K$-maps
by Lemma~\ref{silva}\,(b),
if $\K=\R$ of $\C$).
The mappings $L(\wb{\lambda}_i)\!: L(\wb{G}_i)\to L(G)$
form a cone over~$\cT$ and hence
induce a continuous Lie algebra homomorphism
$\Lambda\!: \g=\dl\, L(\wb{G}_i)\to L(G)$,\vspace{-1mm}
determined by $\Lambda\circ \phi_i=L(\wb{\lambda}_i)$.
To see that
$\Lambda$ is surjective,
let a geometric tangent vector
$v\in L(G)=T_1G$ be given, say $v=[\gamma]$
where $\gamma\!: U\to G$ is a smooth map
on an open $0$-neighbourhood $U\sub \K$,
such that $\gamma(0)=1$.
By Proposition~\ref{mfdc},
after shrinking $U$ we may assume
that $\gamma(U)\sub \wb{\lambda}_i(\wb{G}_i)$
for some $n\in \N$, and that
$\gamma=\wb{\lambda}_i\circ \eta$
for some smooth map $\eta\!: U\to \wb{G}_i$.
Then $\Lambda(\phi_i([\eta]))=L(\wb{\lambda}_i)([\eta])=[\wb{\lambda}_i
\circ \eta]=[\gamma]=v$, as desired.
Because
$\g=\bigcup_{i\in I}\im\, \phi_i$
and $\Lambda\circ \phi_i=L(\wb{\lambda}_i)$,
{\em injectivity\/} of $\Lambda$
follows from the injectivity of
the maps
$L(\wb{\lambda}_i)=T_1(\wb{\lambda}_i)$
established in Theorem~\ref{createmfd}.
By the preceding, $\Lambda$ is an isomorphism
of Lie algebras; as both $\g$ and $L(G)\isom \K^s$
are equipped with the finite topology,
$\Lambda$ also is an isomorphism of topological vector spaces.
Hence $L(G)\isom \g=\dl\, L(\wb{G}_i)$\vspace{-1mm}
naturally. The desired direct limit
properties carry over from $\g$ to $L(G)$.\vspace{1mm}

(b) and (c): Assume that $H$ is a $C^\infty_\K$-Lie group
modelled on a topological $\K$-vector space
and $(f_i)_{i\in I}$
a family
of $C^\infty_\K$-homomorphisms $f_i\!: G_i\to H$
such that $(H,(f_i)_{i\in I})$
is a cone over~$\sys$.
Given $j\in I$, for any $i\geq j$
we then have $f_j=f_i\circ \lambda_{i,j}$
and thus $\ker\,\lambda_{i,j}\sub \ker\, f_j$,
entailing that
$\wb{\bigcup_{i\geq j}\ker\, \lambda_{i,j}}\sub
\ker\, f_j$.
In the situation of (b),
this means that $N_j\sub \ker\, f_j$.
In the situation of~(c),
we assume that $H$ is modelled on a complex
locally convex space.
Then $\ker\, f_j$
is a complex Lie subgroup of $G_j$ (Lemma~\ref{auxil}),
which contains
$\bigcup_{i\geq j}\ker\, \lambda_{i,j}$;
hence also $N_j\sub \ker\, f_j$ in~(c).\vspace{-.3mm}
In any case, we deduce that $f_j=\wb{f}_j\circ q_j$
for a homomorphism $\wb{f}_j\!: \wb{G}_j\to H$,
which is $C^\infty_\K$ because $q_j$ admits smooth
local sections.
Then $((\wb{f}_i)_{i\in I},H)$
is a cone over $\wb{\sys}$ and hence
induces a unique $C^\infty_\K$-homomorphism
$f\! : G\to H$ such that $f\circ \wb{\lambda}_i=\wb{f}_i$,
since $(G,(\wb{\lambda}_i)_{i\in I})=\dl\, \wb{\sys}$\vspace{-.8mm}
in the category of $C^\infty_\K$-Lie groups.
Then $f\circ \lambda_i=f\circ \wb{\lambda}_i\circ q_i
=\wb{f}_i\circ q_i=f_i$
for each $i\in I$, and clearly $f$ is determined
by this property.
Thus $(G,(\lambda_i)_{i\in I})=\dl\, \sys$\vspace{-.8mm}
in the category of $C^\infty_\K$-Lie groups.
In the situation of~(b),
the universal
property of direct limit in the other categories
of interest can be proved by trivial
adaptations of the argument just given.\vspace{1mm}

(d) To establish the first assertion,
we may assume that $I=\N$, and after replacing
$\sys$ by a system equivalent to
$\wb{\sys}$ we may assume
that $G_1\sub G_2\sub \cdots$, each $\lambda_{n,m}$
being the respective inclusion map.
Then $L(G)=\bigcup_{n\in \N}L(G_n)$.
If
$\gamma\!: \R\to L(G)$ is a smooth curve,
then for each $k\in \Z$,
there exists $n_k\in \N$ such that
the relatively compact set $\gamma(]k-1,k+2[)$ is contained in
$L(G_{n_k})$. The finite-dimensional Lie group~$G_{n_k}$
being $c^\infty_\R$-regular, we find a smooth curve
$\eta_k\!: \;]k-1,k+2[\to G_{n_k}$
such that $\eta_k(k)=1$ and
$\delta^r(\eta_k)=\gamma_k$.
We define $\eta(t):=\eta_k(t)\eta_{k-1}(k)\cdots
\eta_1(2)\eta_0(1)$
for $t\in [k,k+1]$ with $k\geq 0$,
and $\eta(t):=\eta_k(t)
\eta_k(k+1)^{-1}\cdots\eta_{-2}(-1)^{-1}
\eta_{-1}(0)^{-1}$
for $t\in [k,k+1]$ with $k<0$.
Then $\eta\!: \R\to G$ is a smooth curve
such that $\eta(0)=1$ and $\delta^r(\eta)=\gamma$.
Thus every $\gamma\in C^\infty(\R,L(G))$
has a right product integral
$\Evol^r_G(\gamma):=\eta\in C^\infty(\R,G)$.
We define
\[
\evol^r_G\!: C^\infty(\R, L(G))\to G,\qquad
\evol^r_G(\gamma):=\Evol^r_G(\gamma)(1)\,.
\]
To see that $\evol^r_G$ is $c^\infty_\R$, let
$\sigma\!: \R\to C^\infty(\R,L(G))$ be a smooth curve.
Given $t_0\in \R$,
let $U\sub \R$ be a
relatively compact neighbourhood
of~$t_0$.
We show that
$\evol^r_G\circ \sigma\!: \R\to G$
is smooth on~$U$.
The evaluation map $C^\infty(\R,L(G))\times \R\to L(G)$,
$(\gamma,t)\mto \gamma(t)$
being continuous (cf.\
Thm.\,3.4.3 and Prop.\,2.6.11 in \cite{Eng}),
$\sigma(U)([-1,2])$
is a compact subset of $L(G)$
and hence contained in $L(G_n)$
for some $n\in \N$, by Lemma~\ref{base}\,(d).
We now consider
\[
\tau\!: U\to C^\infty(]{-1},2[, L(G_n)),\quad \tau(t):=\sigma(t)|_{]{-1},2[}^{L(G_n)}\,.
\]
Then $\tau$ is smooth, because
the restriction map
$C^\infty(\R,L(G))\to C^\infty(]{-1},2[,L(G))$
is continuous linear, and
$C^\infty(]{-1},2[,L(G_n))$
is a closed vector subspace of
$C^\infty(]{-1},2[,L(G))$.
The group $G_n$ being regular,
$\evol^r_{G_n}\!: C^\infty(]{-1},2[,L(G_n))\to G_n$
is smooth.
Since
$\evol^r_G\circ \sigma|_U=\evol^r_{G_n}\circ \tau$
apparently,
we deduce that $\evol^r_G\circ \sigma|_U$
is smooth.
Thus $\evol^r_G$ is $c^\infty_\R$.
The desired direct limit properties can be proved
as in (a) and (b), based on
Proposition~\ref{mfde}.\vspace{1mm}

(e) As a consequence
of Proposition~\ref{mfde},
the $C^\infty_\C$-Lie group~$G$
also is a $c^\infty_\C$-Lie group.
It is $c^\infty_\C$-regular because
its underlying $c^\infty_\R$-Lie group
is $c^\infty_\R$-regular by (d).
The desired direct limit property can be
proved as in (a) and (b).\vspace{1mm}

(f) Using the $c^\omega_\R$-analogue
of Proposition~\ref{mfdf} (see
Proposition~\ref{mfdg}), we see as
in the proof of (a)
that the $C^\infty_\R$-compatible $c^\omega_\R$-manifold structure on~$G$
from Proposition~\ref{mfdg}
turns $G$ into a $c^\omega_\R$-Lie group.
By (d), the latter is $c^\infty_\R$-regular.
To see that it is $c^\omega_\R$-regular,
let $\gamma\!: \R\to L(G)$ be a $c^\omega_\R$-curve
and $\eta:=\Evol^r_G(\gamma)$
be its right product integral.
Using Proposition~\ref{mfdc}
and its $c^\omega_\R$-analogue
(Proposition~\ref{mfdg}),
for each $k\in \N$
we find $n\in \N$ such that $\gamma([-k,k])\sub
L(G_n)$, $\eta([-k,k])\sub G_n$,
and such that $\sigma:=\gamma|_{]-k,k[}^{L(G_n)}$
is $c^\omega_\R$
and $\tau:=\eta|_{]-k,k[}^{G_n}$
smooth.
The finite-dimensional
Lie group~$G_n$
being $c^\omega_\R$-regular,
the product integral $\tau$ of the
$c^\omega_\R$-curve $\sigma$
has to be~$c^\omega_\R$.
Hence $\eta|_{]-k,k[}$
is $c^\omega_\R$ for each $k\in \N$ and thus
$\eta$ is~$c^\omega_\R$.
Hence $G$ is $c^\omega_\R$-regular.
The direct limit property can be established as in~(b).
\end{proof}
We needed to assume local convexity in
Theorem~\ref{limgps}\,(c)
because the proof of the following simple lemma
depends on local convexity.
\begin{la}\label{auxil}
Let $\phi\!: G\to H$
be a $C^\infty_\C$-homomorphism
from
a finite-dimensional complex Lie group
to a complex Lie group modelled on a
locally convex complex
topological vector space.
Then $\,K:=\ker\phi$ is a complex Lie subgroup
of $G$. The same conclusion holds if $H$
is a $c^\omega_\C$-Lie group and $\phi$ a
$c^\omega_\C$-homomorphism.
\end{la}
\begin{proof}
Being a closed subgroup of~$G$,
$K$ is a real Lie subgroup,
with Lie algebra
\[
\ck=\{X\in L(G)\!: \exp_G(\R X)\sub K\}=\{X\in L(G)\!: \phi(\exp_G(\R X))=\{1\}\,\}\, .
\]
Given $X\in \ck$, the map $f\!: \C\to H$, $f(z):=\phi(\exp_G(zX))$
is complex analytic and $f|_\R=1$,
whence $f=1$ by the Identity Theorem.
Hence $\C X\sub \ck$,
whence $\ck$ is a complex Lie subalgebra of $L(G)$.
Therefore $K$ is a complex Lie subgroup \cite[Ch.\,III, \S4.2,
Cor.\,2]{Bou}.\vspace{-3mm}
\end{proof}
\begin{rem}
(a) In the situation of Remark~\ref{preparecx},
the direct system
$((\C^\times)_{n\in \N},(\sigma)_{n\geq m})$
has the direct limit $\R$ in the category of real
Lie groups,
whereas its direct limit in the category of complex
Lie groups is the trivial group.
Hence the conclusions of Theorem~\ref{limgps}\,(a)
become false in general if we replace the injective
quotient system~$\wb{\sys}$
by a non-injective system~$\sys$.
Also note that $\,\dl\, L(\C^\times)=\C\not =\{0\}=L(\{1\})
=L\big(\dl\, \C^\times\big)$ here.

(b) If each $\lambda_{i,j}$ is injective,
then the direct systems $\sys$ and $\wb{\sys}$
are equivalent, whence Theorem~\ref{limgps}\,(a)
remains valid
when $\wb{\sys}$ is replaced with
$\sys$ and all bars are omitted.
\end{rem}
\begin{prop}
Assume that
$((G_i)_{i\in I},(\lambda_{i,j})_{i\geq j})$
is a countable direct
system of
finite-dimensional Lie groups
over $\K\in \{\R,\C\}$
and injective $C^\infty_\K$-homomorphisms,
with direct limit $(G,(\lambda_i)_{i\in I})$.
Then the following holds:
\begin{itemize}
\item[\rm (a)]
$\exp_G=\dl\, \exp_{G_i}\!:
L(G)=\dl\, L(G_i)\to\dl\, G_i=G$\vspace{-.8mm}
is the exponential map of~$G$,
where $(L(G),(L(\lambda_i))_{i\in I})
=\dl\,\big((L(G_i)),(L(\lambda_{i,j}))\big)$.
The map $\exp_G$ is $c^\omega_\K$.
\item[\rm (b)]
The Trotter Product Formula
$\exp_G(x+y)=\lim_{n\to\infty} \big(
\exp_G(\frac{1}{n}x)\exp_G(\frac{1}{n}y)\big)^n$
holds and
$\exp_G([x,y])=\lim_{n\to \infty}
\big(\exp_G(\frac{1}{n}x)\exp_G(\frac{1}{n}y)\exp_G(-\frac{1}{n}x)
\exp_G(-\frac{1}{n}y)\big)^{n^2}$
$($the Commutator Formula$)$,
for all $x,y\in L(G)$.
\item[\rm (c)]
Let
$(H,(\mu_i))=\dl\,\cT$\vspace{-.8mm}
for a direct system
$\cT=((H_i),(\mu_{i,j}))$
of finite-dimensional $\K$-Lie groups
and injective $C^\infty_\K$-homomorphisms,
and assume that $f_i\!: G_i\to H_i$
are $C^\infty_\K$-homomorphisms.
Then $L\Big(\dl\, f_i\Big)\, =\, \dl\, L(f_i)$.\vspace{-.8mm}
Furthermore, every continuous homomorphism
$G\to H$ is $c^\omega_\R$.
\end{itemize}
\end{prop}
\begin{proof}
(a) By Theorem~\ref{limgps}\,(a),
$L(G)=\dl\, L(G_i)$ as a $c^\omega_\K$-manifold.
The family $(\exp_{G_i})_{i\in I}$
of $c^\omega_\K$-maps being compatible with the direct systems
by naturality of~$\exp$, there is a unique
$c^\omega_\K$-map $\exp_G:=\dl\, \exp_{G_i}$
such that $\exp_G\circ L(\lambda_i)=\lambda_i\circ \exp_{G_i}$
for each~$i$. Given $y\in L(G)$,
there are $j\in I$ and $x\in L(G_j)$
such that $y=L(\lambda_j)(x)$.
Then $\xi\!: \R\to G$,\linebreak
$\xi(t):=\exp_G(tx)=
\lambda_j(\exp_{G_j}(ty))$
is a smooth homomorphism
such that $\xi'(0)=$\linebreak
$L(\lambda_j)(\exp_{G_j}'(0).y)=L(\lambda_j)(y)=x$.
Hence $G$ admits an exponential map
(in the sense of \cite[Defn.\,36.8]{KaM}),
and it
is given by $\exp_G$ from above
and hence~$c^\omega_\K$.\vspace{1mm}

(b) Using (a), the assertions directly follow from
the finite-dimensional case.\vspace{1mm}

(c) By Theorem~\ref{limgps}\,(a),
$\alpha:=\dl\, L(f_i)$\vspace{-.8mm} is a continuous $\K$-Lie
algebra homomorphism.
Abbreviate $f:=\dl\, f_i$.\vspace{-.8mm}
From
$\exp_H\circ \alpha=
\big(\dl\, \exp_{G_i}\big)\circ \big(\dl\, L(f_i)\big)
=\dl\,\big(\exp_{H_i}\circ L(f_i)\big)
=\dl\,\big(f_i\circ \exp_{G_i}\big)=f\circ \exp_G$
we deduce that $\alpha=T_0(\exp_H\circ\alpha)
=T_0(f\circ \exp_G)=L(f)$.\vspace{1mm}

Now suppose that
$h\!: G\to H$ is a continuous homomorphism.
We may assume that $I=\N$
and $G_1\sub G_2\sub \ldots$,
$H_1\sub H_2\sub \cdots$.
After replacing $G$ by its identity component~$G_0$, we may assume
that each $G_n$ is connected.
Using Proposition~\ref{mfdc},\vspace{-.8mm} we find $m(n)\in \N$ such that
$h(G_n)\sub H_{m(n)}$,
and such that $h_n:=h|_{G_n}^{H_{m(n)}}$ is continuous
and hence $C^\omega_\R$.
We may assume that $m(1)<m(2)<\cdots$;
after passing to a cofinal subsequence
of the $H_n$'s, without loss of generality $m(n)=n$
for each~$n$. Thus $h=\dl\,h_n$\vspace{-3mm}
is $c^\omega_\R$, by Theorem~\ref{limgps}\,(f).
\end{proof}
\begin{rem}\label{pathoex}
(a) The exponential map of a direct limit
group need not be well-behaved,
as the example $G=\C^\infty \dsemi_{\alpha} \,\R
=\dl\, (\C^n \dsemi \,\R)$\vspace{-.8mm}
with $\alpha\!: \R\to\Aut(\C^\infty)$,
$\alpha(t)((z_k)_{k\in \N}):=(e^{ikt}z_k)_{k\in \N}$
shows. Here $\exp_G$ fails to be injective on
any $0$-neighbourhood,
and the exponential image
$\,\im(\exp_G)$ is not an identity neighbourhood in~$G$
\cite[Example~5.5]{DIR}.\vspace{1mm}

(b) As a consequence of (a), also the exponential
map $\exp_H$ of the complex analytic
Lie group $H:=\C^\infty \dsemi_{\beta} \, \C=\dl\, (\C^n\dsemi\,
\C)$\vspace{-.8mm}
with $\beta(z)((z_k)_{k\in \N}):=(e^{ikz}z_k)_{k\in \N}$
is not injective on any $0$-neighbourhood.
This settles an open problem
by J. Milnor \cite[p.\,31]{MiP}
in the negative,
who asked whether every complex analytic Lie group
is ``of Campbell-Hausdorff type.''
\end{rem}
\section{Integration of locally finite Lie algebras}
A Lie algebra is {\em locally finite\/}
if every finite subset generates a finite-dimensional
subalgebra.
\begin{thm}\label{intthm}
Let $\g$ be a countable-dimensional
locally finite Lie algebra
over $\K\in \{\R,\C\}$.
Then there exists a
$c^\omega_\K$-regular, $c^\omega_\K$-Lie group~$G$, which also is a regular
$C^\infty_\K$-Lie group in Milnor's sense,
such that $L(G)\isom \g$, equipped with the finite topology.
\end{thm}
\begin{proof}
As $\g$ is locally finite and $\dim_\K(\g)\leq \aleph_0$,
there is an ascending sequence
$\g_1\sub \g_2\sub\cdots$
of finite-dimensional subalgebras
of~$\g$,
with union~$\g$.
For each $n\in \N$, let
$G_n$ be a simply connected $\K$-Lie group
with Lie algebra $L(G_n)\isom \g_n$;
fix an isomorphism $\kappa_n\!: L(G_n)\to \g_n$.
If $m\geq n$, then
the Lie algebra homomorphism $\kappa_{m,n}:=\kappa_m^{-1}\circ \kappa_n\!:
L(G_n)\to L(G_m)$
(corresponding to the inclusion map $\g_n\emb \g_m$)
induces a $C^\omega_\K$-homomorphism
$\phi_{m,n}\!: G_n\to G_m$
such that $L(\phi_{m,n})=\kappa_{m,n}$.
Since
$L(\phi_{k,m}\circ \phi_{m,n})=L(\phi_{k,m})\circ L(\phi_{m,n})
=\kappa_{k,m}\circ
\kappa_{m,n}=\kappa_{k,n}=L(\phi_{k,n})$,
we have $\phi_{k,m}\circ\phi_{m,n}=\phi_{k,n}$ for all $k\geq m\geq n$,
whence $((G_n)_{n\in\N},(\phi_{m,n})_{m\geq n})$ is a
direct system of $C^\omega_\K$-Lie groups.
Now take
$G:=\dl \, G_n$\vspace{-.8mm}
in the category of $c^\omega_\K$-Lie groups.
We shall presently show that,
for each $n$,
the normal subgroup $K_n:=\bigcup_{m\geq n}\ker \phi_{m,n}$
of $G_n$ is discrete.
Hence, by Theorem~\ref{limgps}, $G$ is a
$c^\omega_\K$-regular $c^\omega_\K$-Lie group,
$G=\dl \, G_n/K_n$,\vspace{-.8mm}
and
$L(G)=\dl \, L( G_n/K_n)=\dl \, L( G_n)\isom \dl \, \g_n=\g$.\vspace{-.8mm}
For Milnor regularity,
see Theorem~\ref{Milnorreg}.\vspace{1mm}

{\em Each $K_n$ is discrete\/}:
We show that the closure
$N_n:=\wb{K_n}\sub G_n$ is discrete.
The homomorphism $\phi_{m,n}$ has discrete
kernel for all $m,n\in \N$ with $m\geq n$,
because $L(\phi_{m,n})=\kappa_{m,n}$ is
injective. Now $\ker\phi_{m,n}$ being a discrete
normal subgroup of the connected group $G_n$,
it is central. This entails that $N_n\sub Z(G_n)$,
for each~$n$,
whence $(N_n)_0\sub Z(G_n)_0$ is a vector group
being a connected closed subgroup
of a vector group (Lemma~\ref{vectgp}).
It is obvious from the definitions that
$\phi_{m,n}(K_n)\sub K_m$ for all $m\geq n$,
whence
$\phi_{m,n}(N_n)\sub N_m$
and $\phi_{m,n}((N_n)_0)\sub (N_m)_0$.
Being a continuous homomorphism between
vector groups, $\psi_{m,n}:=\phi_{m,n}|_{(N_n)_0}^{(N_m)_0}$
is real linear.\vspace{-.4mm}
Hence, being a real linear map with
discrete kernel, $\psi_{m,n}$ is
injective.
Thus $(N_n)_0=\wb{\bigcup_{m\geq n}\ker \psi_{m,n}}=\{1\}$,
whence~$N_n$ is discrete.
\end{proof}
Here, we used the following fact:
\begin{la}\label{vectgp}
Let $G$ be a simply connected,
finite-dimensional real Lie group.
Then
$Z(G)_0$ is a vector group,
i.e., $Z(G)_0\isom \R^m$ for some $m\in \N_0$.
\end{la}
\begin{proof}
By L\'{e}vi's Theorem,
$L(G)=\Cr\dsemi \cs$
internally,
where $\Cr$ is the radical of~$L(G)$
and $\cs$ a L\'{e}vi
complement
(\cite{Ser},
Part~I, Ch.\,VI, Thm.\,4.1 or
\cite{Bou}, Ch.\,I, \S6.8, Thm.\,5).
Let
$R$ and $S$ be the
analytic subgroups of $G$ corresponding to
$\Cr$ and $\cs$, respectively.
Then $R$ and $S$ are simply connected,
$R$ is a closed normal subgroup of~$G$,
$S$ a closed
subgroup, and $G=R\dsemi S$ internally
\cite[Kor.\,III.3.16]{HaN}.
Now consider the identity component
$Z(G)_0$ of the centre $Z(G)$ of~$G$.
Let $\pi \!: G\to S$ be the projection
onto~$S$, with kernel $R$.
Then $\pi(Z(G)_0)\sub Z(S)_0=\{1\}$,
entailing that $Z(G)_0\sub R$.
Thus $Z(G)_0$ is an analytic subgroup
of the simply connected solvable Lie group~$R$,
whence $Z(G)_0$ is simply connected
\cite[Satz~III.3.31]{HaN}.
Being a simply connected abelian
Lie group, $Z(G)_0$ is a vector group.
\end{proof}
\section{Extension of sections in principal bundles}
We prove a preparatory result
concerning sections in principal bundles,
which will be used later
to discuss closed subgroups and homogeneous
spaces of direct limit groups.
\begin{la}\label{extsec}
Given $\K\in \{\R,\C\}$,
let $\pi_1\!: E_1\to M_1$
be a
$C^\infty_\K$-map between finite-dimensional
$C^\infty_\K$-manifolds
and $\pi_2\!: E_2\to M_2$ be
a finite-dimensional
$G$-principal bundle of class $C^\infty_\K$
whose structure group~$G$
is a finite-dimensional
$C^\infty_\K$-Lie group.
Let $m_1:=\dim_\K(M_1)$ and $m_2:=\dim_\K(M_2)$.
Assume that $\lambda\!: M_1\to M_2$
is an injective $C^\infty_\K$-immersion
and $\Lambda\!: E_1\to E_2$
a $C^\infty_\K$-map such that
$\pi_2\circ \Lambda=\lambda\circ \pi_1$.
Assume that
$\phi_1\!: \Delta^{m_1}_r(\K)\to U_1\sub M_1$
a chart for $M_1$, where
$r>0$,
and $\sigma_1\!: U_1\to E_1$ a $C^\infty_\K$-section
of $\pi_1$. Then, for every $s\in \;]0,r[$,
there exists a chart $\phi_2\!: \Delta^{m_2}_s\to
U_2\sub M_2$
and a $C^\infty_\K$-section
$\sigma_2\!: U_2\to E_2$ of $\pi_2$ such that
$\phi_2(x,0)=\lambda(\phi_1(x))$
for all $x\in \Delta^{m_1}_s$
and $\sigma_2\circ
\lambda|_W
=\Lambda\circ \sigma_1|_W$, where
$W:=\phi_1(\Delta^{m_1}_s)$.
If $\,\K=\R$
and all of $E_1$, $M_1$, $\pi_1$,
the principal bundle $\pi_2$,
$\lambda$, $\Lambda$,
$\phi_1$ and
$\sigma_1$
are real analytic,
then also
$\phi_2$ and $\sigma_2$
can be chosen as real analytic maps.
\end{la}
\begin{proof}
Since $\lambda\circ \phi_1$ is an injective immersion,
there is a chart $\phi_2\!: \Delta^{m_2}_s\to U_2\sub M_2$
such that $\phi_2(x,0)=\lambda(\phi_1(x))$
for all $x\in \Delta^{m_1}_s$
(Lemma~\ref{extpoly}).\vspace{1mm}

{\em The $C^\infty_\K$-case.}
If $\K=\R$, then
$E_2|_{U_2}$ is trivial as a $G$-principal bundle
of class $C^\infty_\R$,
since $U_2\isom \Delta^{m_2}_s$
is paracompact and contractible
(this is a well-known fact,
which can be proved exactly as
\cite[Cor.\,4.2.5]{Hir}).
If $\K=\C$,
then
$E_2|_{U_2}$ is trivial as a $G$-principal bundle
of class $C^\infty_\C$,
since $U_2\isom \Delta^{m_2}_s(\C)$
is a contractible Stein manifold
\cite[Satz~6]{Gr1}.

{\em Real analytic case.}
Since $U_2\isom \Delta^{m_2}_s$
is $\sigma$-compact and contractible,
$E_2|_{U_2}$ is trivial as a topological
$G$-principal bundle
and therefore also as a real analytic
$G$-principal bundle,
by injectivity of $\vartheta$ in
\cite[Teorema~5]{Tog}.

By the preceding,
in either case, we find a
$C^\infty_\K$ (resp., real analytic)
trivialization $\theta\!: E_2|_{U_2}\to U_2\times G$.
Let $\theta_2\!: E_2|_{U_2}\to G$ be the second coordinate
function of~$\theta$.
Define $\sigma_2:=\zeta\circ \phi_2^{-1}\!: U_2\to E_2$,
where $\zeta\!: \Delta^{m_2}_s\to E_2$ is defined via
\[
\zeta(x,y)\; :=\; \theta^{-1}\bigl(\phi_2(x,y),\;
\theta_2((\Lambda \circ \sigma_1 \circ \phi_1)(x))\bigr)\quad
\mbox{for $x\in \Delta^{m_1}_s$, $y\in \Delta^{m_2-m_1}_s$.}
\]
Then $\sigma_2\!: U_2\to E_2$ is a $C^\infty_\K$-section
(resp., $C^\omega_\R$-section) with the required properties. 
\end{proof}
\section{$\!\!$Fundamentals of Lie theory for direct limit groups}
In this section, we develop the basics
of Lie theory for direct limit groups.
Throughout the following,
$G_1\sub G_2\sub \cdots$
is an ascending sequence of finite-dimensional
Lie groups over $\K\in \{\R,\C\}$,
such that the inclusion maps $\lambda_{n,m}\!: G_m\to G_n$
are $C^\omega_\K$-homo\-mor\-phisms,
and $G:=\bigcup_{n\in \N}G_n$,
equipped with the $c^\omega_\K$-Lie group structure
such that $G=\dl\, G_n$\vspace{-.8mm} in the category of
$c^\omega_\K$-Lie groups
(and $C^\infty_\K$-Lie groups).
\begin{numba}
A map $f\!: M\to N$ between $c^\omega_\K$-manifolds
is called {\em $c^\omega_\K$-final\/} if a map $g\!: N\to Z$
into a $c^\omega_\K$-manifold is $c^\omega_\K$
if and only if $g\circ f$ is $c^\omega_\K$.
The map $f$ is {\em $c^\omega_\K$-initial\/}
if a map $g\!: Z\to M$ from a $c^\omega_\K$-manifold~$Z$
to~$M$ is $c^\omega_\K$ if and only if $f\circ g$ is $c^\omega_\K$.
Obvious adaptations
are used to define
$c^\infty_\K$-final, $C^r_\K$-final, $c^\infty_\K$-initial,
and $C^r_\K$-initial maps, where $r\in \N_0\cup\{\infty\}$.
A subset $M\sub N$ of a $c^\omega_\K$-manifold $M$
is called a (split) submanifold
if there exists a (complemented) closed
vector subspace $F$ of the modelling space~$E$
of~$N$ such that, for every $x\in M$, there
exists a chart $\phi\!: U\to V\sub N$ of~$N$ around~$x$
such that $\phi(U\cap F)=M\cap V$. Then $M$, with the induced topology,
can be made a
$c^\omega_\K$-manifold modelled on~$F$,
in an apparent way.
\end{numba}
\begin{prop}[Subgroups are Lie groups]\label{subgplie}
Every subgroup $H\sub G$ admits
a $c^\omega_\K$-Lie group structure
with Lie algebra $L(H)=\{v \in L(G)\!: \exp_G(\K v )\sub H\}=:\ch$
which makes the inclusion map
$\lambda\!: H\to G$ a $C^\infty_\K$-homomorphism
and both a $c^\omega_\K$-initial and a $c^\infty_\K$-initial
map,
and such that $L(\lambda)\!: L(H)\to L(G)$
is an embedding of topological $\K$-vector spaces.
Furthermore, $H=\dl\, H_n$\vspace{-.8mm}
in the category of
$c^\omega_\K$-Lie groups,
where $H_n:=H\cap G_n$ is equipped with
the finite-dimensional $\K$-Lie group structure
induced by $G_n$.
\end{prop}
\begin{proof}
We equip $H_n$ with the finite-dimensional
$C^\omega_\K$-Lie group
structure induced by $G_n$,
which makes the inclusion
map $\lambda_n\!: H_n\to G_n$ an
immersion and a $C^\omega_\K$-initial
and $C^\infty_\K$-initial map
inside the category of finite-dimensional
$C^\omega_\K$- and $C^\infty_\K$-manifolds, respectively 
(see \cite{Bou}, Ch.\,III, \S4.5,
Defn.\,3 and Prop.\,9).
Then the inclusion maps
$\mu_{n,m}\!: H_m\to H_n$
are $C^\omega_\K$-immersions
for $n \geq m$;
we give $H=\bigcup_{n\in \N}H_n$
the $c^\omega_\K$-Lie group structure
such that $(H,(\mu_n)_{n\in \N})
=\dl\,((H_n)_{n\in \N},(\mu_{n,m}))$\vspace{-.8mm}
in the category of $c^\omega_\K$-Lie groups,
where $\mu_n\!: H_n\to H$ is the inclusion map
(see Theorem~\ref{limgps}).
Then $\lambda=\dl\,\lambda_n\!: H\to G$\vspace{-.8mm}
is $c^\omega_\K$.
We have $L(H)=\bigcup_{n\in \N}L(H_n)$
(with obvious identifications)
and $L(G)=\bigcup_{n\in \N}L(G_n)$
by Theorem~\ref{limgps}\,(a),
and each of
$L(\lambda_n)\!: L(H_n)\to L(G_n)$
as well as $L(\lambda)\!: L(H)\to L(G)$
is the respective inclusion map.
Thus $L(\lambda)$ is injective.
Being an injective linear
map between locally convex spaces equipped
with their finest locally
convex topologies, $L(\lambda)$
is a topological embedding
(cf.\ \cite[Prop.\,7.25\,(ii)]{HaM}).
Clearly $L(H)\sub \ch$. If $v\in \ch$, then
$v\in L(G_n)$ for some~$n$ and thus $\exp_{G_n}(\K v)\sub
G_n\cap H=H_n$, whence $v\in L(H_n)\sub L(H)$. Thus $L(H)=\ch$.
Now assume that
$M$ is a $c^\infty_\K$-manifold
and $f\!: M\to H$ a map such that
$\lambda\circ f\!: M\to G$ is
$c^\infty_\K$.
Then, for every smooth map
$\gamma\!: \K\supseteq U\to M$
on an open $0$-neighbourhood $U\sub \K$,
the composition $f\circ \gamma$ maps some $0$-neighbourhood
$V\sub U$ into
some $G_n$ and $(f\circ\gamma)|_V^{G_n}$ is $C^\infty_\K$,
by Proposition~\ref{mfdc}.
Since $f(\gamma(V))\sub G_n\cap H=H_n$ and $H_n$
is $C^\infty_\K$-initial
for maps from finite-dimensional $C^\infty_\K$-manifolds,
we deduce that $(f\circ\gamma)|_V^{H_n}$
is $C^\infty_\K$, whence $(f\circ \gamma)|_V$
is $c^\infty_\K$.
This entails that $f$ is $c^\infty_\K$.
Thus $\lambda\!: H\to G$ is $c^\infty_\K$-initial.
Similarly,
$\lambda$ is $c^\omega_\K$-initial.\vspace{-2mm}
\end{proof}
Note that every $c^\omega_\C$-map between open subsets
of $\C^\infty$ is~$c^\omega_\R$.\,\footnote{In view of
Lemma~\ref{silva}\,(b)
and the $c^\omega_\C$-analogue of
Proposition~\ref{mfdc}
(Proposition~\ref{mfdg}),
this follows from the fact
that holomorphic maps
between open subsets of finite-dimensional
vector spaces are real analytic.}
Hence every direct limit
$G=\dl\,G_n$\vspace{-.8mm}
of complex finite-dimensional Lie groups
has an underlying $c^\omega_\R$-Lie group~$G_\R$.
\begin{la}\label{realcx}
If $\K=\C$
in the situation of Proposition~{\rm \ref{subgplie}},
define $\ch$ as before
and $\cs:=$
$\{v \in L(G)\!: \exp_G(\R v )\sub H\}$.
Let $S_n$ be $H_n$, equipped
with the $C^\omega_\R$-Lie group structure
induced by the $C^\omega_\R$-Lie group
$(G_n)_\R$ underlying $G_n$,
and define $S:=\dl\, S_n$.\vspace{-.8mm}
Thus $S=H$ as a set and an abstract group,
and $\id\!: H_\R\to S$ is $c^\omega_\R$.
Then $\ch=\cs$ $($as a set or real Lie algebra$)$
if and only
if $\,(H_n)_\R=S_n$ $($as a real Lie group$)$ for each $n\in \N$,
if and only if $\,H_\R=S$ $($as a $c^\omega_\R$-Lie group$)$.
\end{la}
\begin{proof}
If $\ch=\cs$,
then for every $n\in \N$ we have
$L(S_n)+iL(S_n)\sub L(H_m)$
for some $m\geq n$. Let $v\in L(S_n)$.
Then $\exp_{G_n}(\C v)=\exp_{G_m}(\C v)
\sub H_m\cap G_n=H_n$,
entailing that $v \in L(H_n)$.
Thus $L(S_n)\sub L(H_n)$ and hence $L(S_n)=L(H_n)$,
whence $S_n=(H_n)_\R$.\vspace{1mm}

If $S_n=(H_n)_\R$ for each $n\in \N$, then
$(\dl\, H_n)_\R=\dl\, (H_n)_\R=\dl\, S_n$,\vspace{-1mm} by
Theorem~\ref{limgps}.\vspace{1mm}

Now suppose that $H_\R=S$. We have
$\ch\sub \cs$ by definition.
If $v\in \cs$, then $\exp_G(\R v)\sub H$
and $\xi\!: \R\to S$, $\xi(t)=\exp_G(tv)=\exp_S(tv)$
is a $c^\omega_\R$-homomorphism.
Since $S=H_\R=\dl\, (H_n)_\R$\vspace{-.8mm} (see Theorem~\ref{limgps}),
Proposition~\ref{mfdc}
entails that $\im(\xi)\sub (H_n)_\R$ for some
$n\in \N$ and that $\xi|^{(H_n)_\R}$
is $C^\omega_\R$. Hence
$\xi=\exp_{H_n}(\sbull w)$
for some $w\in L(H_n)$,
where $w=v$ clearly and thus $\exp_G(\C v)=
\exp_{H_n}(\C v)\sub H_n\sub H$,
whence $v\in \ch$. Therefore $\cs=\ch$.\vspace{-3mm}
\end{proof}
\begin{numba}\label{prepa}
We now specialize
to the case where $H$ is a {\em closed\/} subgroup of $G$;
if $\K=\C$, we assume that
$\cs:=\{v \in L(G)\!:
\exp_G(\R v )\sub H\}$
is a complex Lie subalgebra of~$L(G)$.
Then
the finite-dimensional
$C^\omega_\K$-Lie group structure
induced by~$G_n$ on its closed subgroup
$H_n:=G_n\cap H$
makes $H_n$ a closed
$C^\omega_\K$-submanifold of~$G_n$
(this is obvious in the real
case, and follows for $\K=\C$ using
Lemma~\ref{realcx}).
For each $n\in \N$, we give $G_n/H_n$ the
finite-dimensional $C^\omega_\K$-manifold structure
making the canonical quotient map
$q_n\!: G_n\to G_n/H_n$
a submersion.
Let $q_{n,m}\!: G_m/H_m\to G_n/H_n$
be the uniquely determined $C^\omega_\K$-maps
such that $q_{n,m}\circ q_m=q_n\circ \lambda_{n,m}$.
Then
$\cT:=((G_n/H_n)_{n\in \N},(q_{n,m})_{n\geq m})$
is a direct system
of paracompact, finite-dimensional
$C^\omega_\K$-manifolds
and injective $C^\omega_\K$-immersions,
whence $(M,(\psi_n)_{n\in \N})
:=\dl\,\cT$\vspace{-.8mm}
exists in the category of $c^\omega_\K$-manifolds
(Proposition~\ref{mfdg}).
We have $(M,(\psi_n)_{n\in \N})=\dl\,\cT$\vspace{-.8mm}
also
in the categories of $C^\infty_\K$-manifolds,
and the category of sets.
Consider the quotient map $q\!: G\to G/H$
and the inclusion maps $\lambda_n\!: G_n\to G$.
For each $n\in \N$, the map $q \circ \lambda_n$
factors to an injective map $\mu_n\!: G_n/H_n\to G/H$,
determined by $\mu_n\circ q_n=q\circ\lambda_n$. 
Then $(G/H,(\mu_n)_{n\in \N})$
is a cone over~$\cT$,
and hence induces a map $\mu\!: M\to G/H$.
It is easy to see that $\mu$ is a bijection;
we give $G/H$ the $c^\omega_\K$-manifold structure
making $\mu$ a $c^\omega_\K$-diffeomorphism;
thus $G/H\isom \dl\, G_n/H_n$.\vspace{-.8mm}
Then also $(G/H,(\mu_n)_{n\in \N})=\dl\, \cT$
in the category of $c^\omega_\K$-manifolds.
Since $q=\dl\, q_n$,\vspace{-1mm} the map $q$ is $c^\omega_\K$.
\end{numba}
\begin{prop}[Closed subgroups, quotient groups,
homogeneous spaces]\label{homog}$\;$\\
Let $H$ be a closed subgroup of $G$;
if $\K=\C$, assume that $\{v \in L(G)\!:
\exp_G(\R v)\sub H\}$
is a complex Lie subalgebra of~$G$.
Equip $G/H$ with
the $c^\omega_\K$-manifold structure
described in {\bf \ref{prepa}};
thus $G/H\isom \dl\, G_n/H_n$\vspace{-.8mm}
as a $c^\omega_\K$-manifold.
Then the following holds:
\begin{itemize}
\item[\rm (a)]
$q\!: G\to G/H$ admits local $c^\omega_\K$-sections,
i.e., $q\!: G\to G/H$ is an $H$-principal
bundle of class $c^\omega_\K$.
Therefore $q$ is $c^\omega_\K$-final,
$c^\infty_\K$-final and $C^r_\K$-final, for each $r\in \N_0\cup\{\infty\}$.
\item[\rm (b)]
$H$ is a closed, split $c^\omega_\K$-submanifold
of~$G$.
The $c^\omega_\K$-submanifold structure
on~$H$ and the $c^\omega_\K$-manifold structure
underlying the $c^\omega_\K$-Lie group
structure induced by $G$ on~$H$ $($as described in
Proposition~{\rm \ref{subgplie})\/}
coincide. This manifold structure makes
the inclusion map $H\to G$
a $c^\omega_\K$-initial,
$c^\infty_\K$-initial, and $C^r_\K$-initial map,
for each $r\in \N_0\cup\{\infty\}$.
If $L(H)=\{0\}$, then $H$ is discrete in the topology induced by~$G$.
\item[\rm (c)]
If $H$ is furthermore
a normal
subgroup of $G$,
then the $c^\omega_\K$-manifold structure on
$G/H$
makes the quotient group~$G/H$ a $c^\omega_\K$-regular
$c^\omega_\K$-Lie group
such that
$G/H
=\dl\,G_n/H_n$\vspace{-.8mm}
in the category of $c^\omega_\K$-Lie groups.
\end{itemize}
\end{prop}
\begin{proof}
(a) Let $x\in G/H$;
then there exists $k\in \N$
and $y\in G_k$
such
that $x=q(y)$. Define $z:=q_k(y)$.
Define $r_n:=1+2^{-n}$ for $n\in \N$,
and $d_n:=\dim_\K(G_n/H_n)$.
There exists a $C^\omega_\K$-section
$\tau\!: V\to G_k$ of $q_k$ on some
open neighbourhood $V$ of $z$ in $G_k/H_k$,
and a chart $\phi_k\!: \Delta^{d_k}_{r_k}\to
U_k\sub G_k/H_k$ such that $U_k\sub V$;
we define $\sigma_k:=\tau|_{U_k}$.
Inductively, Lemma~\ref{extsec}
provides charts $\phi_n\!: \Delta^{d_n}_{r_n}\to U_n\sub G_n/H_n$
and $C^\omega_\K$-sections $\sigma_n\!: U_n\to G_n$
such that $q_{n,m}\circ \phi_m|_{\Delta^{d_m}_{r_n}}
=\phi_n|_{\Delta^{d_m}_{r_n}}$
for all $n\geq m\geq k$
and $\sigma_n(q_{n,m}(w))
=\sigma_m(w)$
for all $w\in \phi_m(\Delta^{d_m}_{r_n})$.
Define $W_n:=\phi_n(\Delta^{d_n}_1)$
for $n\in \N$, $n\geq k$.
Then $W:=\bigcup_{n\geq k}\mu_n(W_n)$
is an open subset of~$G/H$, and $(W,(\mu_n|_{W_n}^W)_{n\geq k})=
\dl\,\cW$\vspace{-.8mm} as a $c^\omega_\K$-manifold,
where
$\cW:=((W_n)_{n\geq k},(q_{n,m}|_{W_m}^{W_n}))$.
Now $\sigma:=\dl\, (\sigma|_{W_n})\!: W
=\dl\, W_n\to \dl\, G_n=G$\vspace{-.8mm}
is a $c^\omega_\K$-map,
and it is a section for $q$ because
$q\circ \sigma=\dl\,(q_n\circ \sigma_n|_{W_n})=\dl\, j_n=j$,\vspace{-.8mm}
where $j_n\!: W_n\to G_n/H_n$ and $j\!: W\to G/H$
are the inclusion maps. The remainder is obvious.\vspace{1mm}

(b) For the $c^\omega_\K$-Lie group structure
induced by $G$ on~$H$,
we have
$H=\dl\, H_n$\vspace{-.8mm}
by Proposition~\ref{subgplie},
and this then also holds for
the underlying $c^\omega_\K$ and $c^\infty_\K$-manifold
structures (Theorem~\ref{limgps} (d)--(f)).
Hence,
by the proof of Theorem~\ref{createmfd},
there exists a chart of $H$ around~$1$
of the form $\phi=\dl\,\phi_n\!: P\to Q\sub H$,\vspace{-.8mm}
where, for each~$n$,
the map $\phi_n\!: P_n\to Q_n\sub H_n$ is a chart
of $H_n$ around~$1$, defined
on an open subset $P_n\sub \K^{h_n}$ (where $h_n:=\dim_\K(H_n)$),
$P:=\bigcup_{n\in \N} P_n\sub \K^t$
(where $t:=\sup\{h_n\!: n\in \N\}\in \N_0\!\cup\!\{\infty\}$),
and $Q:=\bigcup_{n\in \N} Q_n\sub H$.
By the proof of Part\,(a) of the present
proposition,
there exist charts $\psi_n\!: \Delta^{d_n}_1\to W_n$
onto open neighbourhoods
$W_n\sub G_n/H_n$ of $q_n(1)$
(where $d_n:=\dim_\K(G_n/H_n)$)
and $C^\omega_\K$-sections $\sigma_n\!: W_n\to G_n$
of~$q_n$, such that $q_{n,m}(W_m)\sub W_n$,
$q_{n,m}\circ \psi_m=\psi_n|_{\Delta^{d_m}_1}$,
and $\sigma_n\circ q_{n,m}|_{W_m}=\sigma_m$
for all $m,n\in \N$ such that $n\geq m$.
Define $V_n:=\im(\sigma_n)Q_n\sub G_n$ and
\[
\theta_n\!: \Delta^{d_n+h_n}_1
\to V_n,\quad
\theta_n(x,y):=\sigma_n(\psi_n(x))\phi_n(y)\quad
\mbox{for $x\in \Delta^{d_n}_1$, $y\in \Delta^{h_n}_1$.}
\]
Since $\sigma_n$ is a section of $q_n$,
the map $\theta_n$ is easily
seen to be injective.
Using the inverse function theorem,
one verifies that $V_n$ is open in~$G_n$
and that $\theta_n$ is a $C^\omega_\K$-diffeomorphism
onto~$V_n$.
Then $V:=\bigcup_{n\in \N}V_n$
is open in $G$, and $\theta:=\dl\,\theta_n\!:
\dl\, \Delta^{d_n+h_n}_1\to V$\vspace{-.8mm}
is a $c^\omega_\K$-diffeomorphism.
Let $\zeta\!: \K^{s+t}=\dl\, \K^{d_n+h_n}\to \K^s\times \K^t$\vspace{-.8mm}
be the natural isomorphism of topological
vector spaces (cf.\ Proposition~\ref{mfdf}),
and $\Omega:=\zeta(\bigcup_n \Delta^{d_n+h_n}_1)\sub \K^s\times
\K^t$.
Then $\kappa:=\theta\circ \zeta^{-1}|_\Omega\!: \Omega\to V$
is a $C^\infty_\K$-diffeomorphism.
By construction of $\theta$,
we have $V\cap H=Q$ and
$\kappa^{-1}(V\cap H)=\Omega\cap (\{0\}\times \K^t)$,
where
$\{0\}\times \K^t$
is a closed,
complemented vector subspace
of $\K^s\times \K^t$.
Hence
$H$ is a split $c^\omega_\K$-submanifold of~$G$.
As the restriction of
$\kappa$ to a submanifold chart of~$H$
is the given chart~$\phi$ of~$H$,
the submanifold structure and the above
manifold structure on~$H$ coincide.
If $L(H)=0$, then the topology on
the Lie group $H$ is discrete and hence
so is the topology on~$H$ as a submanifold
of~$G$, the induced topology.\vspace{1mm}

(c) By construction,
$(G/H,(\mu_n))=\dl\, ((G_n/H_n), (q_{n,m}))$\vspace{-.8mm}
as a $c^\omega_\K$-manifold.
Since each $q_{n,m}$ also is a homomorphism,
Theorem~\ref{limgps} shows that there is a
group structure on the $c^\omega_\K$-manifold
$G/H$ making it a Lie group, and such that each $\mu_n$
becomes a homomorphism. This requirement entails
that $q\!: G\to G/H$ is a homomorphism, whence the
group structure on $G/H$ is the one of the quotient group.
For the second assertion, see Theorem~\ref{limgps}.
\end{proof}
Proposition~\ref{homog}\,(a) entails
that the surjection
$q$ is an open map.
Hence $q$ is a topological quotient map,
and the manifold $G/H$ carries the quotient topology.
\begin{example}
If $G_n$ closed in $G_{n+1}$
for each $n$ and
$K_n\sub G_n$ a maximal compact subgroup
such that $K_1\sub K_2\sub \cdots$,
then
$K:=\bigcup_{n\in \N}K_n$ is a closed subgroup
of~$G=\bigcup_{n\in \N}G_n$.
In fact,
$K_m\cap G_n=K_n$ for $m\geq n$ by maximality,
whence $K\cap G_n=K_n$
is closed in~$G_n$.
\end{example}
\begin{prop}\label{kernels}
If $f\!: G\to A$
is $C^\infty_\K$-
$($resp., $c^\infty_\K$-\,$)$ homomorphism
from $G=\bigcup_{n\in \N}G_n$
into a $C^\infty_\K$-Lie group
modelled on a locally convex space $($resp., a $c^\infty_\K$-
Lie group\,$)$, then
$H:=\ker(f)$ satisfies
the hypotheses of Proposition~{\rm \ref{homog}},
and $L(H)=\ker L(f)$.
\end{prop}
\begin{proof}
In the complex case,
$H\cap G_n=\ker (f|_{G_n})$
is a complex Lie subgroup of~$G_n$
by Lemma~\ref{auxil},
whence the specific hypothesis
of Proposition~\ref{homog}
is satisfied, by
Lemma~\ref{realcx}.
If $w\in \ker L(f)$,
then $\xi\!: \R \to H$, $\xi(t):=f(\exp_G(tw))$
is a $C^\infty_\R$- (resp., $c^\infty_\R$-)
homomorphism
such that $\xi'(0)=L(f)(w)=0$
and thus $\xi=1$
(\cite[La.\,7.1]{Mil}, \cite[La.\,36.7]{KaM}).
Hence $\exp_G(\K w)\sub H$
and therefore $w\in L(H)=\{v\in L(G)\!: \exp_G(\K v)\sub H\}$.
The inclusion $L(H)\sub \ker L(f)$
is trivial.
\end{proof}
\begin{prop}[Canonical Factorization]\label{limmaps}
Let $f\!: G\to A$ be a $c^\infty_\K$-homomorphism
between direct limit groups,
where $G$ is connected,
$G=\bigcup_{n\in \N}G_n$, and
$A=\bigcup_{n\in \N}A_n$.
Equip $G/\ker(f)$
with the $c^\omega_\K$-Lie group structure
from Proposition~{\rm \ref{homog}\,(c)},
and $\im(f)$ with the $c^\omega_\K$-Lie group
structure induced by~$A$ $($as in Proposition~{\rm \ref{subgplie})}.
Let $\wb{f}\!: G/\ker(f)\to \im(f)$
be the bijective homomorphism
induced by~$f$.
Then $\wb{f}$ is a $c^\omega_\K$-diffeomorphism.
\end{prop}
\begin{proof}
In view of Proposition~\ref{mfdd}, we may assume that each $G_n$
and $A_n$ is connected.
Note that $\wb{f}$ is $c^\omega_\K$ because
the inclusion map $\im(f)\to A$ is $c^\omega_\K$-initial
and the quotient map $G\to G/\ker(f)$ is
$c^\omega_\K$-final.
Replacing $G$ with $G/\ker(f)$ and $A$ with $\im(f)$,
we may therefore assume
that $f$ is bijective, and have to
show that $f^{-1}$ is $c^\omega_\K$.
Then $L(f)$ is injective, by Proposition~\ref{kernels}.\vspace{1mm}

{\em $L(f)$ is surjective.}
To see this, let $x \in L(A)=\bigcup_{n\in \N}L(A_n)$;
define $\cs:=\K x$ and $S:=\exp_A(\cs )$.
If $x\not\in \im(L(f))$,
then $\ch_n\cap \cs=\{0\}$
for each $n\in \N$, where
$\ch_n:=L(f)(L(G_n))$.
Given~$n$,
there exists $m\in \N$ such that $L(A_m)\supseteq \ch_n\cup \cs$.
Let
$H_n$ and $S_n$ be the analytic
subgroups of $A_m$
with Lie algebras $\ch_n$
and $\cs$, respectively.
Then $S=S_n$ as a set,
and the group
$H_n\cap S=H_n\cap S_n$ is countable,
because $\ch_n\cap \cs=\{0\}$.
Thus $S=\bigcup_{n\in \N} (S\cap H_n)$
is countable. But $S$ is uncountable, contradiction.
Therefore $x\in \im(L(f))$.\vspace{1mm}

{\em $f^{-1}\!$ is $\,c^\omega_\K$.}
As $A=\dl\, A_n$,\vspace{-1mm}
it suffices to show that $f^{-1}|_{A_n}$
is $c^\omega_\K$, for each $n\in \N$.
Fix~$n$. There exists $m\in \N$
such that $L(f)(L(G_m))\supseteq L(A_n)$.
Let $B$ be the analytic subgroup
of $G_m$ with Lie algebra $L(f)^{-1}(L(A_n))$.
Then $f|_B^{A_n}$ is a bijective
$C^\omega_\K$-homomorphism between
connected finite-dimensional
$\K$-Lie groups
and therefore an isomorphism
of $C^\omega_\K$-Lie groups.
Hence $f^{-1}|_{A_n}^B$ is $C^\omega_\K$
and hence so is $f^{-1}|_{A_n}$.
\end{proof}
\begin{prop}[Universal covering group]\label{coverings}
If $G_n$ is connected for each $n\in \N$,\linebreak
let $\pi_n\!: \wt{G}_n\to G_n$
be the universal covering group,
and $\wt{\lambda}_{n,m}\!: \wt{G}_m\to\wt{G}_n$
be the $C^\infty_\K$-homomorphism
which lifts $\lambda_{n,m}\circ \pi_m$.
Then
$((\wt{G}_n)_{n\in \N},(\wt{\lambda}_{n,m}))$
is a direct system in the category of $C^\infty_\K$-Lie groups;
let
$(\wt{G},(\Lambda_n)_{n\in \N})$ be its direct
limit.
Then $\wt{G}$ is simply connected,
and the $C^\infty_\K$-homomorphism
$\pi:=\dl\, \pi_n\!: \wt{G}\to G$\vspace{-.8mm}
is a universal covering map.
\end{prop}
\begin{proof}\footnote{We cannot use Remark~\ref{homotopy}
because
the $\wt{\lambda}_{n,m}$'s
need not be injective.}
As any connected
$C^\infty_\K$-Lie group,
$G$ has a universal covering
group
$p\!: P\to G$;
thus $G$ is a simply connected
$C^\infty_\K$-Lie group and
$p$ a $C^\infty_\K$-homomorphism
with discrete kernel.
Being a regular topological space and locally diffeomorphic
to $L(G)$, $P$ is smoothly Hausdorff
and hence also is a $c^\infty_\K$-Lie group.
By \cite[Thm.\,38.6]{KaM}, $P$ is
a $c^\infty_\K$-regular Lie group.
Let $\lambda_n\!: G_n\to G$ be the inclusion map.
Since~$P$ is $c^\infty_\K$-regular,
$L(\lambda_n)\!: L(\wt{G}_n)=L(G_n)\to L(G)=L(P)$
integrates to a $c^\infty_\K$-homomorphism
$\alpha_n\!: \wt{G}_n\to P$
(Lemma~\ref{useofreg},
\cite[Thm.\,40.3]{KaM}).
Being a cone, $(P,(\alpha_n))$
induces a $c^\infty_\K$-homomorphism
$\alpha\!: \wt{G}\to P$, determined
by $\alpha\circ \Lambda_n=\alpha_n$.
Recall from Theorem~\ref{limgps}
that
$\wt{G}=\dl\, \wt{G}_n/D_n$,\vspace{-1mm}
where $D_n:=\ker(\Lambda_n)$
and where the limit map
$\mu_n\!: \wt{G}_n/D_n\to \wt{G}$
is obtained by factoring
$\Lambda_n$ over $\wt{G}_n\to \wt{G}_n/D_n$.
Because
$\pi\circ \Lambda_n=\lambda_n \circ \pi_n$,
the subgroup $D_n\sub \ker(\pi_n)$
is discrete.
Hence $L(\wt{G})=\dl\, L(\wt{G}_n/D_n)=\dl\, L(\wt{G}_n)=\dl\,
L(G_n)=L(G)$.\vspace{-1mm}
It is easily verified that
$L(\alpha)=\id_{L(G)}$ with respect to these
identifications.
Now
$\wt{G}$
being $c^\infty_\K$-regular and~$P$
simply connected, $\id\!: L(P)=L(G)\to L(G)=L(\wt{G})$
induces a $c^\infty_\K$-homomorphism
$\beta\!: P\to \wt{G}$,
determined by $L(\beta)=\id_{L(G)}$.
Since $L(\alpha\circ \beta)=\id_{L(G)}=L(\id_P)$,\vspace{-.3mm}
we have $\alpha\circ \beta=\id_P$ by
\cite[La.\,7.1]{Mil}.
Likewise, $\beta\circ \alpha=\id_{\wt{G}}$.
Thus $\wt{G}\isom P$ is the universal covering group.
\end{proof}
\begin{prop}[Integration of Lie algebra homomorphisms]\label{inthom}
Assume that  $G=\bigcup_{n\in \N}G_n$ is simply connected.
Then the following holds:
\begin{itemize}
\item[\rm (a)]
Every $\K$-Lie algebra
homomorphism $\alpha\!: L(G)\to L(H)$
into the Lie algebra of
a $c^\infty_\K$-regular
$c^\infty_\K$-Lie group~$H$
integrates to a $c^\infty_\K$-homomorphism
$\beta\!: G\to H$ such that $L(\beta)=\alpha$.
If $\K=\R$ and $H$ is a $c^\omega_\R$-regular $c^\omega_\R$-Lie group
here, then $\beta$ is $c^\omega_\R$.
\item[\rm (b)]
Every $\K$-Lie algebra
homomorphism $\alpha\!: L(G)\to L(H)$
into the Lie algebra of a $\K$-analytic BCH-Lie group
$($see {\rm \cite{GCX})} integrates to a $C^\infty_\K$-
$($and $c^\omega_\K$-$)$ homomorphism
$\beta\!: G\to H$.
\end{itemize}
\end{prop}
\begin{proof}
(a) See \cite[Thm.\,40.3]{KaM} and Lemma~\ref{useofreg}.

(b) Let $((\wt{G}_n),(\wt{\lambda}_{n,m}))$,
$\,(\wt{G}, (\Lambda_n))$,
$\,\pi_n\!: \wt{G}_n\to G_n$,
and $\pi\!: \wt{G}\to G$
be as in Proposition~\ref{coverings}.
Because $G$ is simply connected,
the covering homomorphism~$\pi$
is an isomorphism. Hence $\wt{G}=G$ and
$\Lambda_n=\lambda_n\circ \pi_n$
without loss of generality, where
$\lambda_n\!: G_n\to G$ is the inclusion map.
Now $H$ and $\wt{G}_n$
being BCH, the homomorphism $\alpha_n:=\alpha \circ L(\Lambda_n)$
integrates to a $C^\omega_\K$-homomorphism $\beta_n\!:
\wt{G}_n\to H$~\cite[Prop.\,2.8]{GCX}.
Then $(H,(\beta_n))$ is a cone
and hence induces a $C^\infty_\K$-
(and $c^\omega_\K$-) homomorphism
$\beta\!: G\to H$ such that
$\beta\circ \Lambda_n=\beta_n$.
Clearly $L(\beta)=\alpha$.
\end{proof}
\begin{prop}[Integration of Lie subalgebras]\label{intsubalg}
\,Given a $\K$-Lie subalgebra $\ch$
of\linebreak
$L(G)$,
equip the subgroup $H:=\langle \exp_G(\ch)\rangle$
with the $c^\omega_\K$-Lie group structure
described in Proposition~{\/\rm \ref{subgplie}}.
Then $H$ is connected, and $L(H)=\ch$.
Furthermore,
$H=\dl\, H_n$\vspace{-.8mm}
where $H_n:=\langle \exp_{G_n}(\ch_n)\rangle$
is the analytic subgroup of~$G_n$
with Lie algebra $\ch_n:=\ch \cap L(G_n)$.
\end{prop}
\begin{proof}
Consider the inclusion map
$f\!: S\to G$, where
$S:=\dl\, H_n=\bigcup_{n\in \N}H_n$.\vspace{-.8mm}
Then $S=H$ as an abstract group.
We have $L(S)=\bigcup_{n\in \N}L(H_n)=\ch$,
and $f$ is $c^\omega_\K$ because
each $f|_{H_n}$ is so.
By Proposition~\ref{limmaps},
$f$ is a $c^\omega_\K$-diffeomorphism onto
$\im(f)=H$,
equipped with
the $c^\omega_\K$-Lie group structure
induced by~$G$. Thus~$S=H$
as $c^\omega_\K$-Lie groups.
\end{proof}
Before we can discuss universal
complexifications of direct limit groups,
we need to re-examine universal
complexifications of finite-dimensional Lie groups.
\begin{la}\label{help}
Let $G$ be a finite-dimensional real Lie group,
and $\gamma_G\!: G\to G_\C$
be its universal complexification in the category
of finite-dimensional complex Lie groups.
Let $\alpha \!: G\to H$ be a $c^\infty_\R$-homomorphism
from~$G$ to a $c^\omega_\C$-regular $c^\omega_\C$-Lie group~$H$.
Then there exists a unique
$c^\omega_\C$-homomorphism
$\beta\!: G_\C\to H$ such that $\beta\circ \gamma_G=\alpha$.
\end{la}
\begin{proof}
{\em We assume first that $G$ is connected.}
Let $p\!: \wt{G}\to G$ be the universal
covering group of~$G$ and $S$ be a simply connected
complex Lie group with Lie algebra $L(G)_\C$.
Let $\lambda\!: L(G)\to L(G)_\C$ be the inclusion map
and $\kappa\!: \wt{G}\to S$ be the unique
$c^\infty_\R$-homomorphism such that $L(\kappa)=\lambda$.
Set $\Pi:=\ker(p)\isom \pi_1(G)$
and let $N\sub S$ be the smallest closed
complex Lie subgroup such that $\kappa(\Pi)\sub N$.
Let $q\!: S\to S/N=:G_\C$ be the canonical quotient map.
Then there exists a $c^\infty_\R$-homomorphism
$\gamma_G\!: G\to G_\C$
such that $\gamma_G\circ p=q\circ \kappa$.\vspace{1mm}

Let $\alpha\!: G\to H$ be a $c^\infty_\R$-homomorphism
into a $c^\omega_\C$-regular
$c^\omega_\C$-Lie group~$H$.
Then Lemma~\ref{useofreg}
provides a $c^\omega_\C$-homomorphism
$\eta\!: S\to H$ such that $L(\eta)$
is the $\C$-linear extension of~$L(\alpha)$.
Then $\eta\circ \kappa=\alpha\circ p$,
because $L(\eta\circ \kappa)=L(\alpha)=L(\alpha\circ p)$.
Thus $\kappa(\Pi)\sub \ker(\eta)$,
where $\ker(\eta)$ is a closed, complex Lie subgroup
of~$S$ by Lemma~\ref{auxil}.
Thus $N\sub \ker(\eta)$, and thus $\eta$ factors
to a $c^\omega_\C$-homomorphism $\beta\!: G_\C=S/N\to H$
such that $\beta\circ q=\eta$.
From $\beta\circ \gamma_G\circ p=\beta\circ q\circ \kappa=
\eta\circ \kappa=\alpha\circ p$
we deduce that $\beta\circ\gamma_G=\alpha$,
and clearly $\beta$ is uniquely determined
by this property. By the preceding, $\gamma_G\!: G\to G_\C$
is a universal complexification of the
$c^\infty_\R$-Lie group~$G$ in the category
of $c^\omega_\C$-regular
$c^\omega_\C$-Lie groups; since $G_\C$ is finite-dimensional,
$\gamma_G\!: G\to G_\C$
also is the universal complexification of~$G$
in the category of finite-dimensional complex Lie groups.\vspace{1mm}

{\em If $G$ is not necessarily connected\/},
then its identity component $G_0$ has a
universal complexification in
the category $\cA$ of $c^\omega_\C$-regular
$c^\omega_\C$-Lie groups,
which is finite-dimensional.
As in \cite[Prop.\,5.2]{GCX},
we see that the $c^\infty_\R$-Lie group
$G$ has a universal complexification
$\gamma_G\!: G\to G_\C$
in~$\cA$, and $(G_\C)_0$
is a universal complexification for
$G_0$ and therefore finite-dimensional.
Hence $G_\C$ is finite-dimensional,
and hence it coincides with the universal
complexification of~$G$ in the category
of finite-dimensional complex Lie groups.
\end{proof}
\begin{prop}[Universal complexifications]\label{univcx}
Let $\sys:=((G_n)_{n\in \N},(\lambda_{n,m})_{n\geq m})$
be a direct system of finite-dimensional real Lie groups
and $C^\omega_\R$-homomorphisms,
$(G,(\lambda_n)):=\dl\,\sys$\vspace{-.8mm}
in the category of $c^\omega_\R$-Lie groups
and
$(G_\C,(\kappa_n)_{n\in \N}):=\dl\, \big(((G_n)_\C)_{n\in \N},
((\lambda_{n,m})_\C)\big)$
in the category of $c^\omega_\C$-Lie groups,
where
$\gamma_n\!:G_n\to (G_n)_\C$
is a universal complexification
for $G_n$ in the category of finite-dimensional
complex Lie groups,
and $(\lambda_{n,m})_\C\!: (G_m)_\C\to (G_n)_\C$
the uniquely determined complex analytic homomorphism
such that $(\lambda_{n,m})_\C\circ \gamma_m=\gamma_n\circ
\lambda_{n,m}$.
Let $\gamma_G:=\dl\,\gamma_n\!: G\to G_\C$.\vspace{-1.6mm}
Then the following holds:
\begin{itemize}
\item[\rm (a)]
$\gamma_G\!: G\to G_\C$
is a universal complexification
of the $c^\infty_\R$-Lie group~$G$
in the category of $c^\omega_\C$-regular
$c^\omega_\C$-Lie groups
in the sense that
for every $c^\infty_\R$-homomorphism
$\alpha\!: G\to H$ into a $c^\omega_\C$-regular $c^\omega_\C$-Lie group
$H$, there exists a uniquely determined
$c^\omega_\C$-homomorphism $\beta\!: G_\C\to H$
such that $\beta\circ \gamma_G=\alpha$.
\item[\rm (b)]
$\gamma_G|_{G_0}^{(G_\C)_0}$
is the universal complexification of $G_0$,
and the map $G/G_0\to G_\C/(G_\C)_0$,
$xG_0\mto \gamma_G(x)(G_\C)_0$ is a bijection.
\item[\rm (c)]
If $G$ is simply connected, then $\gamma_G$
has discrete kernel.
\item[\rm (d)]
If $\gamma_G$ has discrete kernel,
then $L(G_\C)=L(G)_\C$,
$\,\im\,\gamma_G$
is closed in~$G_\C$,
and $\gamma_G|^{\smim\, \gamma_G}$
is a local $c^\omega_\R$-diffeomorphism
onto $\im\,\gamma_G$,
equipped with the $c^\omega_\R$-Lie group
structure induced by~$(G_\C)_\R$.
\end{itemize}
\end{prop}
\begin{proof}
(a) By Lemma~\ref{help},
for each $n\in \N$ there exists a unique
$c^\omega_\C$-homomorphism $\beta_n\!: (G_n)_\C\to H$
such that $\beta_n\circ \gamma_n=\alpha\circ \lambda_n$.
Clearly $(H,(\beta_n))$ is a cone,
whence there exists a unique $c^\omega_\C$-homomorphism
$\beta\!: G_\C =\dl\, (G_n)_\C\to H$\vspace{-.8mm}
such that $\beta\circ \kappa_n=\beta_n$.
Then $\beta\circ \gamma_G=\alpha$,
and it is easily verified that $\beta$ is uniquely determined
by this property.\vspace{1.5mm}

(b) Compare \cite[Prop.\,5.2]{GCX}.\vspace{1.5mm}

(c) Using
Theorem~\ref{intthm}, we find a simply connected,
$c^\omega_\C$-regular
$c^\omega_\C$-Lie group $S$ with Lie algebra $L(S)=L(G)_\C$.
As $G$ is simply connected,
the inclusion map $j\!: L(G)\emb L(G)_\C$
integrates to a $c^\infty_\R$-homomorphism
$\eta\!: G\to S$ such that $L(\eta)=j$.
Then $\eta\!: G\to S$
is a universal complexification for~$G$
in the category of $c^\omega_\C$-regular $c^\omega_\C$-Lie
groups (cf.\ proof of Lemma~\ref{help}
or \cite[La.\,IV.4]{GaN}).
Let $K:=\ker(\eta)=\ker(\gamma_G)$.
Because $L(\eta)=j$ is injective,
we have $L(K)=\ker L(\eta)=\{0\}$
(Proposition~\ref{kernels}).
Hence $K$ is discrete
when equipped with the
real Lie group structure induced by~$G$.
The topology on the latter
coincides with the
topology
induced by~$G$, as~$K$ is closed
(Prop.\,\ref{homog}\,(b)).
Hence $\ker(\gamma_G)=K$ is discrete.\vspace{2mm}

(d) Since $\ker(\gamma_G)$ is discrete,
$L(\gamma_G)$ is injective (Proposition~\ref{kernels}),
enabling us to identify $L(G)$ with $\im\, L(\gamma_G)$
as a real locally convex space.
Let $(G_\C)_\ops$ be $G_\C$, equipped with the opposite
complex structure; by the universal property of $G_\C$,
there is a unique $c^\omega_\C$-homomorphism $\sigma\!:G_\C\to
(G_\C)_\ops$ such that $\sigma\circ \gamma_G=\gamma_G$.
We now consider $\sigma$ as an antiholomorphic
self-map of $G_\C$.
Thus $L(\sigma)$ is $\C$-antilinear.
As in \cite[La.\,IV.2]{GaN}, we see
that $\sigma$ is an involution.
We have $L(G)\sub L(G_\C)^\sigma$
for the fixed space of $L(\sigma)$.
Since $L(G_\C)=$\linebreak
$L(G)+iL(G)$ by construction
of $G_\C$, it easily follows that
$L(G_\C)=L(G)\oplus i L(G)=L(G)_\C$
and thus $L(G)=L(G_\C)^\sigma$.
We now give the closed subgroup $(G_\C)^\sigma:=\Fix(\sigma)$
the $c^\omega_\R$-Lie group structure induced by $(G_\C)_\R$.
Then $\gamma_G(G)\sub (G_\C)^\sigma$,
and it is easy to see that $L((G_\C)^\sigma):=\{
v\in L(G_\C)\!: \exp_{G_\C}(\R v)\sub (G_\C)^\sigma\}=L(G)$.
Thus $C:=((G_\C)^\sigma)_0=\langle\exp_{G_\C}(L(G))\rangle
=\gamma_G(G_0)$,
and now Proposition~\ref{limmaps}
entails that $\gamma_G|_{G_0}^C$
is a local $c^\omega_\R$-diffeomorphism.
To complete the proof,
note that $(G_\C)_0\cap \gamma_G(G)=\gamma_G(G_0)=C$
by~(b),
whence $\gamma_G(G)$ is a locally closed subgroup
of~$G_\C$
and hence closed.\vspace{-1.5mm}
\end{proof}
\section{Proof of regularity in Milnor's sense}
\begin{thm}\label{Milnorreg}
Every direct limit group $G=\dl\,G_n$\vspace{-1mm}
over $\K\in \{\R,\C\}$
is a regular $C^\infty_\R$-Lie group in Milnor's sense.
More precisely,
for every $k\in \N\cup\{\infty\}$,
every $C^k_\R$-curve
$\gamma\!:
[0,1]\to G$ admits a right product integral
$\eta=\Evol^r_G(\gamma)\in C^{k+1}([0,1],G)$
such that $\eta(0)=1$,
and the corresponding
right evolution map
\[
\evol^r_G\!: C^k([0,1], L(G))\to G,\quad \evol^r_G(\gamma):=
\Evol^r_G(\gamma)(1)
\]
is $C^\infty_\K$ and $c^\omega_\K$.
\end{thm}
\begin{proof}
Fix $k$. The strategy of the proof
is as follows.
First, we show that product integrals exist
and that $\evol^r_G$ is continuous.
Next, we
show that $\evol^r_G$ is complex analytic
if $\K=\C$.
Finally,
for $\K=\R$,
we deduce smoothness of $\evol^r_G$
from the smoothness of $\evol^r_{G_\C}$.\vspace{1mm}

{\bf Step~1.}
Since $\evol^r_G$ takes its values in the connected
component of~$G$,
we may assume that $G$ is connected.
Using that
$\delta^r(p\circ \gamma)=\delta^r\gamma$
for curves in $\wt{G}$
(cf.\ \cite[38.4\,(3)]{KaM}),
where $p\!: \wt{G}\to G$
is the universal
covering map,
we may assume that $G$ is simply connected.
Furthermore, we may assume that $G=\bigcup_{n\in \N}G_n$,
where $G_1\sub G_2\sub \cdots$
and each $G_n$ is connected.
Let $j_n\!: G_n\to G$
be the inclusion map.
We abbreviate $d_n:=\dim_\K(G_n)$, $s:=\sup\{d_n\!: n\in \N\}$
and let $\phi=\dl\, \phi_n\!: P\to Q$\vspace{-.8mm}
be a chart of $G$
around~$1$,
where $P=\bigcup_{n\in \N}\Delta^{d_n}_2$,
$Q:=\bigcup_{n\in \N} Q_n$ and
$\phi_n\!: \Delta^{d_n}_2\to Q_n$ is a chart of~$G_n$
around~$1$,
such that $\phi_n(0)=1$.
We identify $L(G_n)=T_1(G_n)$ with $\K^{d_n}$
using the chart $\phi_n$,
and $L(G)$ with $\K^s$ using~$\phi$;
then $L(j_n)\!: \K^{d_n}\to \K^s$
is the inclusion map, for each $n\in \N$.\vspace{1mm}

{\bf Step~2:} {\em $\evol^r_G$ exists.} To see this, let
$\gamma\in C^k([0,1],L(G))$.
Then there exists $n\in \N$
such that $\im\, \gamma\sub L(G_n)$.
Then $\gamma|^{L(G_n)}$
is~$C^k$.
It is a standard fact
(based on the local existence and uniqueness
of solutions to differential equations)
that there exists
$\eta\in C^{k+1}([0,1],G_n)$
such that $\delta^r\eta=\gamma|^{L(G_n)}$.
Then $\Evol^r_G(\gamma):=j_n\circ \eta$ is~$C^{k+1}$
and $\delta^r(j_n\circ \eta)=L(j_n)\circ \gamma|^{L(G_n)}=\gamma$.
Thus $\evol^r_G(\gamma)$ exists,
and $\evol^r_G\circ C^k ([0,1],L(j_n))=j_n\circ
\evol^r_{G_n}$.\vspace{1mm}

{\bf Step~3.}
The inclusion map
$C^k([0,1], L(G)) \to C^1([0,1], L(G))$
being continuous linear for each $k$,
it suffices to prove that
$\evol^r_G\!: C^1([0,1],L(G))\to G$
is $C^\infty_\K$ and $c^\omega_\K$.
We may therefore assume that $k=1$
for the rest of the proof.\vspace{1mm}

{\bf Step~4:} {\em $\evol^r_G$ is continuous at nice $\gamma_0$'s.}
We show that $\evol^r_G$
is continuous at
$\gamma_0 \in C^1([0,1],L(G))$,
provided that $\im(\gamma_0)\sub \K^{d_1}=L(G_1)$
and $\im (\eta_0)\sub \phi_1(\Delta^{d_1}_{1/2})$,
where $\eta_0:=\Evol^r_G(\gamma_0)$.
To this end, let
$W$ be an open neighbourhood
of $\evol^r_G(\gamma_0)=\eta_0(1)$ in~$G$;
abbreviate $\zeta_0:=\phi^{-1}_1\circ \eta_0$.
Then $\phi^{-1}(W)$
is an open neighbourhood of $\zeta_0(1)$,
whence $\phi^{-1}(W)-\zeta_0(1)\supseteq \Delta^{d_1}_{\ve_1}\oplus
\bigoplus_{n\geq 2} \Delta^{d_n-d_{n-1}}_{\ve_n}$
for certain $\ve_n>0$;
we may assume that
$1\geq \ve_1\geq \ve_2\geq \cdots$.
Define $r_n:=1-2^{-n}$ for $n\in \N$.
Equip each $\K^{d_n}$ with the supremum norm.
There is $R>0$ such that $\|\gamma_0\|_\infty\leq R$.

For $n\in \N$,
consider the map $f_n\!: \K^{d_n}\times \Delta^{d_n}_2
\to \K^{d_n}$,
$f_n(y,x):=  \frac{d}{ds}\big|_{s=0}\phi_n^{-1}(\phi_n(sy)\phi_n(x))$,
which expresses the
map $L(G_n) \times G_n\to TG_n$, $(y,x)\mto T_1(\rho_x).y$
(with right translation $\rho_x\!: G_n\to G_n$)
in local coordinates (forgetting the fibre).
Then $\zeta_0'(t)=f_n(\gamma_0(t),\zeta_0(t))$
for all $t\in [0,1]$, because $\delta^r(\eta_0)=\gamma_0$.
By compactness of $\wb{\Delta}^{d_n}_1$
and $\wb{\Delta}^{d_n}_{R+n-1}$, there exists $k_n>0$
such that for the operator norms of the partial
differentials we have
\[
\|d_2f_n(v,x,\sbull)\|\leq k_n\quad
\mbox{for all $v\in \Delta^{d_n}_{R+n-1}$ and $x\in \Delta^{d_n}_1$,}
\]
and such that for the operator norms
of the continuous linear maps $f_n(\sbull,x)$
we have $\|f_n(\sbull,x)\|\leq k_n$ for all $x\in \Delta^{d_n}_1$.
Choose $\alpha_n>0$ so small that
\begin{equation}\label{ve}
\frac{\alpha_n}{k_n}\, \bigl(e^{k_n}-1\bigr)\;\leq \; 2^{-n-1}\,\ve_n\,.
\end{equation}
Define $s_n:=\min\{\frac{\alpha_n}{k_n},1\}$.
Suppose that
$\gamma\!: [0,1]\to \Delta^{d_n}_{R+n-1}$
is a $C^1$-curve for which there exists
a $C^1$-curve $\eta\!: [0,1]\to \Delta^{d_n}_{1-2^{-n}}$
solving the initial value problem
$\eta(0)=0$, $\eta'(t)=f_n(\gamma(t),\eta(t))$.
Then
$\|d_2f_n(\gamma(t),x,\sbull)\|\leq k_n$
for all $t\in [0,1]$ and $x\in \Delta^{d_n}_1$.
Let $\wb{\gamma}\!: [0,1]\to \K^{d_n}$
be a $C^1$-curve such that $\|\wb{\gamma}-\gamma\|_\infty< s_n$.
Then $\im (\wb{\gamma})\sub \Delta^{d_n}_{R+n}$,
and
\[
\|f_n(\wb{\gamma}(t),x)-f_n(\gamma(t),x)\|=\|f_n(\wb{\gamma}(t)-\gamma(t),x)\|
\leq \|f_n(\sbull,x)\|\cdot \|\wb{\gamma}(t)-\gamma(t)\|
\leq k_n s_n\leq \alpha_n
\]
for all $x\in \Delta^{d_n}_1$.
Furthermore,
$\eta(t)+y\in \Delta^{d_n}_{1-2^{-n-1}}\sub \Delta^{d_n}_1$
for all $t\in [0,1]$ and $y\in \K^{d_n}$ such that
$\|y\|\leq \frac{\alpha_n}{k_n}(e^{k_n}-1)\leq 2^{-n-1}\ve_n\leq
2^{-n-1}$.
Using
\cite[(10.5.6)]{Die}, we therefore
find a solution $\xi \!: [0,1]\to \Delta^{d_n}_1$
to the initial value problem $\xi(0)=0$,
$\,\xi'(t)=f_n(\wb{\gamma}(t),\xi(t))$,
such that
\begin{equation}\label{contribut}
\|\xi-\eta\|_\infty\; \leq \; \frac{\alpha_n}{k_n}(e^{k_n}-1)\; \leq
\; 2^{-n-1}\ve_n\,.
\end{equation}
Hence $\im (\xi)\sub \Delta^{d_n}_{1-2^{-n-1}}$
in particular.

\noindent
We now define
$\Omega:=
\Delta^{d_n}_{s_1}\oplus\bigoplus_{n\geq 2}\Delta^{d_n-d_{n-1}}_{s_n}$,
considering $\K^s$ as the locally convex direct sum
$\K^{d_1}\oplus \bigoplus_{n\geq 2} \K^{d_n-d_{n-1}}$.
Then
$\gamma_0+C^1([0,1],\Omega)$
is an open neighbourhood of $\gamma_0$ in $C^1([0,1],L(G))$.
Let $\gamma\in \gamma_0+C^1([0,1],\Omega)$.
Then $\gamma-\gamma_0=\sum_{n=1}^\infty\gamma_n$,
where $\gamma_n$ is the coordinate function
taking its values in $\Delta^{d_1}_{s_1}$,
resp., in $\Delta^{d_n-d_{n-1}}_{s_n}$.
There exists $\ell\in \N$ such that $\gamma_n=0$ for all
$n\geq \ell$.
Considering $\gamma_0$, $\gamma_0+\gamma_1$, $\,\ldots\,$,
$\sum_{n=0}^\ell\gamma_n =\gamma$ in turn,
from the existence of $\zeta_0$
we inductively deduce by the preceding arguments
that there exists a solution
$\zeta_n\!: [0,1]\to \Delta^{d_n}_{1-2^{-n-1}}$ to
the initial value problem
$\zeta'_n(t)=f_n\big(\gamma_0(t)+\cdots+\gamma_n(t),\,\zeta_n(t)\big)$,
$\zeta_n(0)=0$,
for $n=1,\ldots, \ell$,
such that $\|\zeta_n-\zeta_{n-1}\|_\infty\leq 2^{-n-1}\ve_n$
(see (\ref{contribut})).
Then $\eta:=\phi\circ \zeta_\ell$
is the right product integral
for $\gamma$, and thus $\evol^r_G(\gamma)=\eta(1)\in W$
because{\small
\[
\phi^{-1}(\eta(1))-\phi^{-1}(\eta_0(1))
=\zeta_\ell(1)-\zeta_0(1)=\!\sum_{n=1}^\ell
(\zeta_n(1)-\zeta_{n-1}(1))
\in \Delta^{d_1}_{\ve_1}\oplus
\bigoplus_{n=2}^\ell \Delta^{d_n-d_{n-1}}_{\ve_n}\!\sub \phi^{-1}(W)\,.
\]}\!
Hence $\evol^r_G$ is indeed continuous at $\gamma_0$.\vspace{1mm}

{\bf Step~5.} {\em $\evol^r_G$ is continuous.}
Let $\wb{\gamma}\in C^1([0,1],L(G))$.
After passing to a subsequence,
we may assume that $\im(\wb{\gamma})\sub \K^{d_1}=L(G_1)$.
Let $\wb{\eta}:=\Evol^r_G(\wb{\gamma})$.
We find a partition $0=t_0<t_1<\ldots <t_\ell=1$
such that
$\wb{\eta}_j([0,1])\sub \phi_1(\Delta^{d_1}_{1/2})$
for each $j\in \{0,\ldots,\ell-1\}$,
where $\wb{\eta}_j\!:[0,1]\to G$,
$\wb{\eta}_j(t)=\wb{\eta}(t_j +t(t_{j+1}-t_j))\, \wb{\eta}(t_j)^{-1}$.
Then $\wb{\eta}_j=\Evol^r_G(\wb{\gamma}_j)$,
where $\wb{\gamma}_j\!: [0,1]\to L(G)$,
$\wb{\gamma}_j(t):=(t_{j+1}-t_j)\cdot
\wb{\gamma}(t_j+t(t_{j+1}-t_j))$ are mappings
which satisfy the hypotheses
of Step~4. Thus $\evol^r_G$ is continuous
at $\wb{\gamma}_j$.
Since $\gamma\mto \gamma_j$
is continuous and
$\evol^r_G(\gamma)
=\evol^r_G(\gamma_{\ell-1})\cdots \evol^r_G(\gamma_1)
\evol^r_G(\gamma_0)$,
we deduce that $\evol^r_G$ is continuous at~$\wb{\gamma}$.\,\footnote{We have
even established continuity with respect to the topology of
uniform convergence\,!}\vspace{1mm}

{\bf Step~6:\/}
{\em $\evol^r_G$ is $C^\infty_\C$ if $\K=\C$.}
It suffices to show that $\evol^r_G$
is $C^\infty_\C$ on some open neighbourhood
of each $\gamma_0\in C^1([0,1], L(G))$
such that $\gamma_0([0,1])\sub L(G_1)$
and such that $\eta_0:=\Evol^r_G(\gamma_0)$
has image in $\phi_1(\Delta^{d_1}_{1/2})$,
by arguments similar to those just employed.
Let $\Omega$ be as in Step~4,
and $U:=\gamma_0+C^1([0,1],\Omega)$.
As shown in Step~4,
$\eta:=\Evol^r_G(\gamma)$
has image in
$Q=\im(\phi)$, for each $\gamma\in U$,
$\zeta:=\phi^{-1}\circ \eta$
satisfies $\zeta(0)=0$, and
$\zeta'(t)=f_n(\gamma(t),\zeta(t))$
for each $n$ such that $\zeta([0,1])\sub \C^{d_n}$.
Now suppose that $\gamma\in U$
and $\theta\in C^1([0,1],L(G))$.
There exists $n$ (which we fix now)
such that $\gamma,\theta$
have image in~$\C^{d_n}$.
Then $\sigma_z:=\gamma+z\theta\in U$ for $z$
in some $0$-neighbourhood $V\sub \C$,
and $\im(\sigma_z)\sub \C^{d_n}$
for each $z\in V$. Let $\tau_z:=\phi^{-1}\circ
\Evol^r_G(\sigma_z)$.
Then $\tau_z$ solves the
initial value problem $\tau_z(0)=0$, $\,\tau_z'(t)=f_n(\sigma_z(t),\tau_z(t))$.
Consider
$f\!: [0,1]\times \Delta^{d_n}_1 \times V\to \C^{d_n}$,
$f(t,x,z):=f_n(\sigma_z(t),x)$.
Then $f(t,x,z)=f_n(\gamma(t),x)+zf_n(\theta(t),x)$,
showing that the differentiability requirements
of
\cite[Thm.\,3.6.1]{Car}
are satisfied.\footnote{To apply
the theorem, note that $f$
extends to an open set,
because $\gamma$ and $\theta$
extend to open intervals by Borel's theorem.\\[7.8cm]}
Hence $u(t,z):=\tau_z(t)$ is $C^1_\R$
in $(t,z)$ on an open neighbourhood
of $I\times \{0\}$ in $I\times V$,
and the map $h\!: [0,1]\to \cL_\R(\C,\C^{d_n})$,
$h(t):=d_2u(t,0,\sbull)$
to the space of $\R$-linear maps $\C\to \C^{d_n}$
is $C^1_\R$ and solves the initial
value problem
\begin{equation}\label{recx}
h(0)=0, \quad h'(t)=b(t)\circ h(t)+c(t),
\end{equation}
where
$c(t)(z)=z\cdot f_n(\theta(t),\tau_0(t))$
and $b(t)=d_2f_n(\sigma_0(t),\tau_0(t),\sbull)$.
Since $b(t)\in \cL_\C(\C^{d_n},\C^{d_n})$ actually for each~$t$
and $c(t)\in \cL_\C(\C,\C^{d_n})$,
we can interpret (\ref{recx})
also as a linear differential equation
for $\cL_\C(\C,\C^{d_n})$-valued functions.
This implies that $h(t)\in \cL_\C(\C,\C^{d_n})$
for each~$t$, i.e., $h(t)=d_2u(t,0,\sbull)$
is complex linear. Hence
$\frac{d}{dz}\big|_{z=0}\phi^{-1}(\evol^r_G(\gamma+z\theta))
=\frac{d}{dz}\big|_{z=0}\tau_z(1)$
$=\frac{\partial}{\partial z}
\big|_{z=0}u(1,z)$
exists as a complex derivative.\vspace{1mm}

By the preceding, $\psi:=\phi^{-1}\circ \evol^r_G|_U\!: U\to \C^s$
admits complex directional derivatives at each
point. Hence $\psi$ is G-analytic in the sense
of \cite[Defn.\,5.5]{BaS},
by
\cite[Prop.\,5.5]{BaS} and
\cite[Thm.\,3.1]{BaS}.
Being G-analytic and continuous,
$\psi$ is complex analytic \cite[Thm.\,6.1\,(i)]{BaS}.\vspace{1mm}

{\bf Step~7:\/}
{\em $\evol^r_G$ is $C^\infty_\R$ and $c^\omega_\R$ if $\K=\R$.}
Because $G$ is assumed
simply connected, we know that
$H:=\gamma_G(G)$
is a closed subgroup of~$G_\C$,
that
$\gamma_G$ has discrete kernel,
and that $\gamma_G$ is a local
$c^\omega_\R$-diffeomorphism
onto $H$, equipped with
the real Lie group structure induced by
$(G_\C)_\R$ (see Proposition~\ref{univcx}\,(c)
and (d)).
Since $H$ is $C^\infty_\R$-initial in $G_\C$
and $c^\omega_\R$-initial
(Proposition~\ref{homog}\,(b)),
we deduce from the smoothness
(and $c^\omega_\R$-property)
of $\gamma_G\circ \evol^r_G=\evol^r_{G_\C}\circ L(\gamma_G)$
that $\gamma_G|^H\circ\evol^r_G$ is $C^\infty_\R$
and $c^\omega_\R$.
As $\evol^r_G$ is continuous and
$\gamma_G|^H$ a local $C^\infty_\R$- (and $c^\omega_\R$-)
diffeomorphism,
this implies that $\evol^r_G$ is $C^\infty_\R$
and $c^\omega_\R$.
\end{proof}
\noindent
{\footnotesize {\bf Helge Gl\"{o}ckner}, TU~Darmstadt, FB~Mathematik~AG~5,
Schlossgartenstr.\,7, 64289 Darmstadt, Germany.\\
E-Mail: gloeckner@mathematik.tu-darmstadt.de}
\end{document}